%
\newif\ifloadreferences\loadreferencestrue
%
%
%
%
%
\let\myfrac=\frac%
\input eplain %
\let\frac=\myfrac%
\let\myfootnote=\footnote%
\input amstex \input epsf %
\let\footnote=\myfootnote%
%
%
\loadeufm\loadmsam\loadmsbm\message{symbol names}\UseAMSsymbols\message{,}%
%
\font\myfontdefault=cmr10%
\newif\ifmakebiblio%
\newif\ifinappendices%
\newif\ifundefinedreferences%
\newif\ifchangedreferences%
\makebibliofalse%
\undefinedreferencesfalse%
\changedreferencesfalse%
%
%
%
%
%
\def\setcatcodes{\catcode`\!=0 \catcode`\\=11}%
{\global\let\noe=\noexpand%
\catcode`\@=11 \catcode`\_=11 \setcatcodes%
!global!def!_@@internal@@makeref#1{%
!global!expandafter!def!csname #1ref!endcsname##1{%
!csname _@#1@##1!endcsname%
!expandafter!ifx!csname _@#1@##1!endcsname!relax%
    !write16{#1 ##1 not defined - run saving references}%
    !undefinedreferencestrue%
!fi}}%
!global!def!_@@internal@@makelabel#1{%
!global!expandafter!def!csname #1label!endcsname##1{%
!edef!temptoken{!csname #1info!endcsname}%
!ifloadreferences%
!expandafter!ifx!csname _@#1@##1!endcsname!relax%
!write16{#1 ##1 not hitherto defined - rerun saving references}%
!changedreferencestrue%
!else%
!expandafter!ifx!csname _@#1@##1!endcsname!temptoken%
!else%
!write16{#1 ##1 reference has changed - rerun saving references}%
!changedreferencestrue%
!fi%
!fi%
!else%
!expandafter!edef!csname _@#1@##1!endcsname{!temptoken}%
!edef!textoutput{!write!references{\global\def\_@#1@##1{!temptoken}}}%
!textoutput%
!fi}}%
!global!def!makecounter#1{!_@@internal@@makelabel{#1}!_@@internal@@makeref{#1}}%
!unsetcatcodes%
}
%
%
%
%
%
\def\turnintolatin#1{\ifcase #1 _\or i\or ii\or iii\or iv\or v\or vi\or vii\or viii\or ix\or x\or xi\or xii\or xiii\or xiv\or xv\or xvi\or xvii\or xviii\or xix\or xx\or xxi\or xxii\or xxiii\or xxiv\or xxv\or xxvi\fi}%
\def\alphanum#1{\ifcase #1 _\or A\or B\or C\or D\or E\or F\or G\or H\or I\or J\or K\or L\or M\or N\or O\or P\or Q\or R\or S\or T\or U\or V\or W\or X\or Y\or Z\fi}%
\newwrite\references%
\ifloadreferences{\catcode`\@=11 \catcode`\_=11 \global\def\_@citation@Ahlfors{1}
\global\def\_@citation@AndEtAl{2}
\global\def\_@citation@BarbFill{3}
\global\def\_@citation@Bers{4}
\global\def\_@citation@BonMonSch{5}
\global\def\_@citation@BonSeppi{6}
\global\def\_@citation@Buser{7}
\global\def\_@citation@FillSmi{8}
\global\def\_@citation@FillSmiII{9}
\global\def\_@citation@GilbTrud{10}
\global\def\_@citation@Jost{11}
\global\def\_@citation@LabI{12}
\global\def\_@citation@LabII{13}
\global\def\_@citation@Mess{14}
\global\def\_@citation@Mumford{15}
\global\def\_@citation@NewNir{16}
\global\def\_@citation@Schlenker{17}
\global\def\_@citation@Schoen{18}
\global\def\_@citation@SchoenYau{19}
\global\def\_@citation@Tamburelli{20}
\global\def\_@citation@Thurston{21}
\global\def\_@citation@TrapaniValli{22}
\global\def\_@citation@Tromba{23}
\global\def\_@citation@Wolf{24}
\global\def\_@head@Introduction{1}
\global\def\_@subhead@EquivariantIsometricEmbeddings{1.1}
\global\def\_@eqn@MinkowskiMetric{\relax \unhbox \voidb@x \hbox {{\relax \tenrm (1.1)}}}
\global\def\_@eqn@SemiDirectGroupLaw{\relax \unhbox \voidb@x \hbox {{\relax \tenrm (1.2)}}}
\global\def\_@eqn@Equivariance{\relax \unhbox \voidb@x \hbox {{\relax \tenrm (1.3)}}}
\global\def\_@proc@MainTheoremI{1.1.1}
\global\def\_@rmk@MainTheoremI{\relax \unhbox \voidb@x \hbox {1.1.1}}
\global\def\_@subhead@GHMCMinkowskiSpacetimes{1.2}
\global\def\_@subhead@OneHarmonicMaps{1.3}
\global\def\_@eqn@TVEnergyFunctional{\relax \unhbox \voidb@x \hbox {{\relax \tenrm (1.4)}}}
\global\def\_@proc@TrapaniValliI{1.3.1}
\global\def\_@eqn@CodazziFieldCondition{\relax \unhbox \voidb@x \hbox {{\relax \tenrm (1.5)}}}
\global\def\_@proc@TrapaniValliII{1.3.2}
\global\def\_@proc@IntroCompactness{1.3.3}
\global\def\_@rmk@IntroCompactness{\relax \unhbox \voidb@x \hbox {1.3.1}}
\global\def\_@rmk@IntroCompactnessII{\relax \unhbox \voidb@x \hbox {1.3.2}}
\global\def\_@subhead@Acknowledgements{1.4}
\global\def\_@head@Preliminaries{2}
\global\def\_@subhead@LinearAlgebraicPreliminaries{2.1}
\global\def\_@eqn@ActionOfGL{\relax \unhbox \voidb@x \hbox {{\relax \tenrm (2.1)}}}
\global\def\_@eqn@BundleProjection{\relax \unhbox \voidb@x \hbox {{\relax \tenrm (2.2)}}}
\global\def\_@eqn@MatrixInnerProduct{\relax \unhbox \voidb@x \hbox {{\relax \tenrm (2.3)}}}
\global\def\_@eqn@DefnComplexStructure{\relax \unhbox \voidb@x \hbox {{\relax \tenrm (2.4)}}}
\global\def\_@eqn@DefnOneZeroSeminorm{\relax \unhbox \voidb@x \hbox {{\relax \tenrm (2.5)}}}
\global\def\_@eqn@JLinAndAntilin{\relax \unhbox \voidb@x \hbox {{\relax \tenrm (2.6)}}}
\global\def\_@eqn@FormulaForComplexLinearComponent{\relax \unhbox \voidb@x \hbox {{\relax \tenrm (2.7)}}}
\global\def\_@eqn@FormulaForSigma{\relax \unhbox \voidb@x \hbox {{\relax \tenrm (2.8)}}}
\global\def\_@eqn@InvarianceOfSigma{\relax \unhbox \voidb@x \hbox {{\relax \tenrm (2.9)}}}
\global\def\_@eqn@TwoDimMatrixRelations{\relax \unhbox \voidb@x \hbox {{\relax \tenrm (2.10)}}}
\global\def\_@eqn@FirstDerivativeA{\relax \unhbox \voidb@x \hbox {{\relax \tenrm (2.11)}}}
\global\def\_@eqn@FirstDerivativeB{\relax \unhbox \voidb@x \hbox {{\relax \tenrm (2.12)}}}
\global\def\_@proc@FirstDerivative{2.1.1}
\global\def\_@subhead@GradientAndDivergence{2.2}
\global\def\_@eqn@DefnGradientFunction{\relax \unhbox \voidb@x \hbox {{\relax \tenrm (2.13)}}}
\global\def\_@eqn@DefnDivergenceVectorField{\relax \unhbox \voidb@x \hbox {{\relax \tenrm (2.14)}}}
\global\def\_@eqn@DefnDivergenceEndmorphism{\relax \unhbox \voidb@x \hbox {{\relax \tenrm (2.15)}}}
\global\def\_@eqn@ObjectsInTermsOfExtAlg{\relax \unhbox \voidb@x \hbox {{\relax \tenrm (2.16)}}}
\global\def\_@eqn@FundamentalRelation{\relax \unhbox \voidb@x \hbox {{\relax \tenrm (2.17)}}}
\global\def\_@proc@FundamentalRelation{2.2.1}
\global\def\_@rmk@FundamentalRelation{\relax \unhbox \voidb@x \hbox {2.2.1}}
\global\def\_@head@PerturbationTheory{3}
\global\def\_@subhead@TheEnergyDensity{3.1}
\global\def\_@eqn@DefnOneZeroEnergyDensity{\relax \unhbox \voidb@x \hbox {{\relax \tenrm (3.1)}}}
\global\def\_@eqn@DefnTotalOneZeroEnergy{\relax \unhbox \voidb@x \hbox {{\relax \tenrm (3.2)}}}
\global\def\_@eqn@ComparisonWithTV{\relax \unhbox \voidb@x \hbox {{\relax \tenrm (3.3)}}}
\global\def\_@subhead@OrbitsOfTheDiffeomorphismGroup{3.2}
\global\def\_@eqn@DefnEnergyOverOrbit{\relax \unhbox \voidb@x \hbox {{\relax \tenrm (3.4)}}}
\global\def\_@eqn@DefnLTwoGradient{\relax \unhbox \voidb@x \hbox {{\relax \tenrm (3.5)}}}
\global\def\_@eqn@FormulaForLieDerivative{\relax \unhbox \voidb@x \hbox {{\relax \tenrm (3.6)}}}
\global\def\_@eqn@DefnDeltaXA{\relax \unhbox \voidb@x \hbox {{\relax \tenrm (3.7)}}}
\global\def\_@eqn@VariationOfEInTermsOfDeltaA{\relax \unhbox \voidb@x \hbox {{\relax \tenrm (3.8)}}}
\global\def\_@eqn@DerivativeOfEnergy{\relax \unhbox \voidb@x \hbox {{\relax \tenrm (3.9)}}}
\global\def\_@eqn@DerivativeOfEnergyAtSymm{\relax \unhbox \voidb@x \hbox {{\relax \tenrm (3.10)}}}
\global\def\_@proc@DerivativeOfEnergyAtSym{3.2.2}
\global\def\_@eqn@CriticalPointCondition{\relax \unhbox \voidb@x \hbox {{\relax \tenrm (3.11)}}}
\global\def\_@proc@CorDerivativeOfEnergyAtSym{3.2.3}
\global\def\_@eqn@SecondVariationFormula{\relax \unhbox \voidb@x \hbox {{\relax \tenrm (3.12)}}}
\global\def\_@proc@SecondVariationFormula{3.2.4}
\global\def\_@subhead@RecoveringEllipticity{3.3}
\global\def\_@eqn@FirstComponentOfOperatorWRTX{\relax \unhbox \voidb@x \hbox {{\relax \tenrm (3.13)}}}
\global\def\_@eqn@DifferentialCurvatureFormulaA{\relax \unhbox \voidb@x \hbox {{\relax \tenrm (3.14)}}}
\global\def\_@eqn@DifferentialCurvatureFormulaB{\relax \unhbox \voidb@x \hbox {{\relax \tenrm (3.15)}}}
\global\def\_@proc@DifferentialCurvatureFormula{3.3.1}
\global\def\_@eqn@DefinitionOfF{\relax \unhbox \voidb@x \hbox {{\relax \tenrm (3.16)}}}
\global\def\_@eqn@DefinitionOfG{\relax \unhbox \voidb@x \hbox {{\relax \tenrm (3.17)}}}
\global\def\_@eqn@ExplicitFormulaForF{\relax \unhbox \voidb@x \hbox {{\relax \tenrm (3.18)}}}
\global\def\_@eqn@DiffInvarianceOfF{\relax \unhbox \voidb@x \hbox {{\relax \tenrm (3.19)}}}
\global\def\_@eqn@DerivativeOfG{\relax \unhbox \voidb@x \hbox {{\relax \tenrm (3.20)}}}
\global\def\_@proc@DerivativeOfG{3.3.2}
\global\def\_@eqn@EllipticOperatorWRTX{\relax \unhbox \voidb@x \hbox {{\relax \tenrm (3.21)}}}
\global\def\_@proc@EllipticityOfOperatorWRTX{3.3.3}
\global\def\_@rmk@EllipticityOfOperatorWRTX{\relax \unhbox \voidb@x \hbox {3.3.1}}
\global\def\_@subhead@PerturbingOneHarmonicDiffeomorphisms{3.4}
\global\def\_@eqn@ControlOfModifiedFunctional{\relax \unhbox \voidb@x \hbox {{\relax \tenrm (3.22)}}}
\global\def\_@proc@ControlOfModifiedFunctional{3.4.1}
\global\def\_@eqn@PerturbationOperator{\relax \unhbox \voidb@x \hbox {{\relax \tenrm (3.23)}}}
\global\def\_@proc@EquivalenceOfModifiedFunctional{3.4.2}
\global\def\_@proc@PerturbationOfOneHarmonicDiffeomorphisms{3.4.3}
\global\def\_@head@Compactness{4}
\global\def\_@subhead@ControllingTheConformalClass{4.1}
\global\def\_@eqn@DefinitionOfIntermediateComplexStructure{\relax \unhbox \voidb@x \hbox {{\relax \tenrm (4.1)}}}
\global\def\_@eqn@DefinitionOfIntermediateMetric{\relax \unhbox \voidb@x \hbox {{\relax \tenrm (4.2)}}}
\global\def\_@eqn@DefinitionOfSystole{\relax \unhbox \voidb@x \hbox {{\relax \tenrm (4.3)}}}
\global\def\_@proc@CompactnessOfIntermediateConformalClasses{4.1.1}
\global\def\_@proc@PrecompactnessForHyperbolicMetrics{4.1.2}
\global\def\_@eqn@ExtremalLengthFormula{\relax \unhbox \voidb@x \hbox {{\relax \tenrm (4.4)}}}
\global\def\_@eqn@DefinitionOfModuleOfConformalStructure{\relax \unhbox \voidb@x \hbox {{\relax \tenrm (4.5)}}}
\global\def\_@eqn@UpperBoundOfConformalModule{\relax \unhbox \voidb@x \hbox {{\relax \tenrm (4.6)}}}
\global\def\_@proc@UpperBoundOfConformalModule{4.1.3}
\global\def\_@rmk@UpperBoundOfConformalModule{\relax \unhbox \voidb@x \hbox {4.1.1}}
\global\def\_@eqn@LowerBoundOfConformalModule{\relax \unhbox \voidb@x \hbox {{\relax \tenrm (4.7)}}}
\global\def\_@proc@LowerBoundOfConformalModule{4.1.4}
\global\def\_@eqn@UncertaintyPrinciple{\relax \unhbox \voidb@x \hbox {{\relax \tenrm (4.8)}}}
\global\def\_@proc@PrecompactnessForComplexStructures{4.1.5}
\global\def\_@subhead@HarmonicMaps{4.2}
\global\def\_@eqn@DefinitionOfLaplacian{\relax \unhbox \voidb@x \hbox {{\relax \tenrm (4.9)}}}
\global\def\_@eqn@LaplacianInCoordinateCharts{\relax \unhbox \voidb@x \hbox {{\relax \tenrm (4.10)}}}
\global\def\_@eqn@DefinitionOfHarmonicity{\relax \unhbox \voidb@x \hbox {{\relax \tenrm (4.11)}}}
\global\def\_@eqn@FHarmonicity{\relax \unhbox \voidb@x \hbox {{\relax \tenrm (4.12)}}}
\global\def\_@eqn@FHarmonicityInCoordinates{\relax \unhbox \voidb@x \hbox {{\relax \tenrm (4.13)}}}
\global\def\_@eqn@DefinitionOfEnergyForFunctions{\relax \unhbox \voidb@x \hbox {{\relax \tenrm (4.14)}}}
\global\def\_@eqn@DefinitionOfEnergyForFunctionsII{\relax \unhbox \voidb@x \hbox {{\relax \tenrm (4.15)}}}
\global\def\_@eqn@CompactnessForFHarmonicFunctions{\relax \unhbox \voidb@x \hbox {{\relax \tenrm (4.16)}}}
\global\def\_@proc@CompactnessForFHarmonicFunctions{4.2.1}
\global\def\_@rmk@CompactnessForFHarmonicFunctionsI{\relax \unhbox \voidb@x \hbox {4.2.1}}
\global\def\_@eqn@PreMaximumPrinciple{\relax \unhbox \voidb@x \hbox {{\relax \tenrm (4.17)}}}
\global\def\_@proc@LiouvilleTheorem{4.2.2}
\global\def\_@eqn@COneCompactness{\relax \unhbox \voidb@x \hbox {{\relax \tenrm (4.18)}}}
\global\def\_@proc@COneCompactness{4.2.3}
\global\def\_@rmk@COneCompactness{\relax \unhbox \voidb@x \hbox {4.2.2}}
\global\def\_@subhead@Compactness{4.3}
\global\def\_@eqn@CompactnessDefinitionOfAm{\relax \unhbox \voidb@x \hbox {{\relax \tenrm (4.19)}}}
\global\def\_@eqn@CompactnessBigPhiAndLittlePhi{\relax \unhbox \voidb@x \hbox {{\relax \tenrm (4.20)}}}
\global\def\_@eqn@CompactnessDefnOfLittlePhi{\relax \unhbox \voidb@x \hbox {{\relax \tenrm (4.21)}}}
\global\def\_@proc@CompactnessTrapaniValli{4.3.1}
\global\def\_@proc@CompactnessForCodazziFields{4.3.2}
\global\def\_@eqn@CompactnessIntermediateConformalStructure{\relax \unhbox \voidb@x \hbox {{\relax \tenrm (4.22)}}}
\global\def\_@eqn@CompactnessDefinitionOfTildeJ{\relax \unhbox \voidb@x \hbox {{\relax \tenrm (4.23)}}}
\global\def\_@eqn@CompactnessDefinitionOfAlpha{\relax \unhbox \voidb@x \hbox {{\relax \tenrm (4.24)}}}
\global\def\_@proc@CompactnessAlphaHarmonicity{4.3.3}
\global\def\_@eqn@CompactnessDefinitionOfPsi{\relax \unhbox \voidb@x \hbox {{\relax \tenrm (4.25)}}}
\global\def\_@eqn@CompactnessDefinitionOfBetaI{\relax \unhbox \voidb@x \hbox {{\relax \tenrm (4.26)}}}
\global\def\_@eqn@CompactnessDefinitionOfBetaII{\relax \unhbox \voidb@x \hbox {{\relax \tenrm (4.27)}}}
\global\def\_@eqn@CompactnessBetaHarmonicityI{\relax \unhbox \voidb@x \hbox {{\relax \tenrm (4.28)}}}
\global\def\_@eqn@CompactnessBetaHarmonicityII{\relax \unhbox \voidb@x \hbox {{\relax \tenrm (4.29)}}}
\global\def\_@eqn@CompactnessEnergyBoundsI{\relax \unhbox \voidb@x \hbox {{\relax \tenrm (4.30)}}}
\global\def\_@eqn@CompactnessEnergyBoundsII{\relax \unhbox \voidb@x \hbox {{\relax \tenrm (4.31)}}}
\global\def\_@proc@MainExistenceThmTP{4.3.4}
\global\def\_@head@TheWeylProblem{5}
\global\def\_@subhead@TheDerivativesOfE{5.1}
\global\def\_@eqn@TwoArgumentEnergyFunctional{\relax \unhbox \voidb@x \hbox {{\relax \tenrm (5.1)}}}
\global\def\_@eqn@SymmetryOfTwoArgumentEnergyFunctional{\relax \unhbox \voidb@x \hbox {{\relax \tenrm (5.2)}}}
\global\def\_@eqn@DefinitionOfFunctionOverTeichmuellerSpace{\relax \unhbox \voidb@x \hbox {{\relax \tenrm (5.3)}}}
\global\def\_@eqn@DefinitionOfGeodesic{\relax \unhbox \voidb@x \hbox {{\relax \tenrm (5.4)}}}
\global\def\_@eqn@FormulaForPhiZero{\relax \unhbox \voidb@x \hbox {{\relax \tenrm (5.5)}}}
\global\def\_@proc@WeylPeterssonGeodesic{5.1.1}
\global\def\_@eqn@DefinitionOfBt{\relax \unhbox \voidb@x \hbox {{\relax \tenrm (5.6)}}}
\global\def\_@eqn@DefinitionOfgt{\relax \unhbox \voidb@x \hbox {{\relax \tenrm (5.7)}}}
\global\def\_@eqn@FamilyOfHt{\relax \unhbox \voidb@x \hbox {{\relax \tenrm (5.8)}}}
\global\def\_@eqn@BasicRelation{\relax \unhbox \voidb@x \hbox {{\relax \tenrm (5.9)}}}
\global\def\_@proc@BasicRelation{5.1.2}
\global\def\_@eqn@FirstDerivativeOfF{\relax \unhbox \voidb@x \hbox {{\relax \tenrm (5.10)}}}
\global\def\_@eqn@FirstDerivativeOfFAtZero{\relax \unhbox \voidb@x \hbox {{\relax \tenrm (5.11)}}}
\global\def\_@proc@FirstDerivativeOfF{5.1.3}
\global\def\_@eqn@UsefulRelationIII{\relax \unhbox \voidb@x \hbox {{\relax \tenrm (5.12)}}}
\global\def\_@proc@UsefulRelationIII{5.1.4}
\global\def\_@rmk@UsefulRelationIII{\relax \unhbox \voidb@x \hbox {5.1.1}}
\global\def\_@eqn@ConvexityCodazziRelation{\relax \unhbox \voidb@x \hbox {{\relax \tenrm (5.13)}}}
\global\def\_@eqn@SecondDerivativeOfF{\relax \unhbox \voidb@x \hbox {{\relax \tenrm (5.14)}}}
\global\def\_@proc@SecondDerivativeOfF{5.1.5}
\global\def\_@proc@StrictConvexityOfE{5.1.6}
\global\def\_@subhead@ProofOfMainTheorem{5.2}
\global\def\_@eqn@SumOfConvexFunctionals{\relax \unhbox \voidb@x \hbox {{\relax \tenrm (5.15)}}}
\global\def\_@eqn@SumOfAOneAndATwoIsExact{\relax \unhbox \voidb@x \hbox {{\relax \tenrm (5.16)}}}
\global\def\_@eqn@PrimitiveOfSumOfAOneAndATwo{\relax \unhbox \voidb@x \hbox {{\relax \tenrm (5.17)}}}
\global\def\_@eqn@DefinitionOfXOneAndXTwo{\relax \unhbox \voidb@x \hbox {{\relax \tenrm (5.18)}}}
\global\def\_@eqn@CharacteristicPropertyOfU{\relax \unhbox \voidb@x \hbox {{\relax \tenrm (5.19)}}}
\global\def\_@eqn@DefinitionOfRhoOneAndRhoTwo{\relax \unhbox \voidb@x \hbox {{\relax \tenrm (5.20)}}}
\global\def\_@eqn@DefinitionOfSupportFunctions{\relax \unhbox \voidb@x \hbox {{\relax \tenrm (5.21)}}}
\global\def\_@eqn@LimitsOfPhisCoincide{\relax \unhbox \voidb@x \hbox {{\relax \tenrm (5.22)}}}
\global\def\_@head@TheGeneralisedLorentzianMetric{A}
\global\def\_@eqn@MinkowskiBilinearForm{\relax \unhbox \voidb@x \hbox {{\relax \tenrm (A.1)}}}
\global\def\_@eqn@LorentzMetricInQuotient{\relax \unhbox \voidb@x \hbox {{\relax \tenrm (A.2)}}}
\global\def\_@eqn@MetricOverQuotient{\relax \unhbox \voidb@x \hbox {{\relax \tenrm (A.3)}}}
\global\def\_@eqn@ConformallyRescaledMetric{\relax \unhbox \voidb@x \hbox {{\relax \tenrm (A.4)}}}
\global\def\_@eqn@CovDerOfRescaledMetric{\relax \unhbox \voidb@x \hbox {{\relax \tenrm (A.5)}}}
\global\def\_@eqn@GeodesicsOfRescaledMetric{\relax \unhbox \voidb@x \hbox {{\relax \tenrm (A.6)}}}
\global\def\_@eqn@DefnOfGeneralisedLorentzMetric{\relax \unhbox \voidb@x \hbox {{\relax \tenrm (A.7)}}}
\global\def\_@eqn@GeodesicOfLorentzianMetric{\relax \unhbox \voidb@x \hbox {{\relax \tenrm (A.8)}}}
\global\def\_@eqn@ExponentialMap{\relax \unhbox \voidb@x \hbox {{\relax \tenrm (A.9)}}}
\global\def\_@eqn@DomainOfExponentialMap{\relax \unhbox \voidb@x \hbox {{\relax \tenrm (A.10)}}}
\global\def\_@eqn@DefinitionOfATilde{\relax \unhbox \voidb@x \hbox {{\relax \tenrm (A.11)}}}
\global\def\_@eqn@DefinitionOfTildeU{\relax \unhbox \voidb@x \hbox {{\relax \tenrm (A.12)}}}
\global\def\_@eqn@DefinitionOfTildePsi{\relax \unhbox \voidb@x \hbox {{\relax \tenrm (A.13)}}}
\global\def\_@head@Bibliography{B}
 }%
\else{\openout\references=references.tex }%
\fi%
%
%
\newcount\headno%
\global\headno=0%
\def\headinfo{\ifinappendices\alphanum\headno\else\the\headno\fi}%
\def\nextheadno{\global\advance\headno by 1 \global\subheadno=0 \global\eqnno=0 \headinfo}%
\makecounter{head}%
%
%
\newcount\subheadno%
\global\subheadno=0%
\def\subheadinfo{\headinfo.\the\subheadno}%
\def\nextsubheadno{\global\advance\subheadno by 1 \global\procno=0 \global\rmkno=0 \subheadinfo}%
\makecounter{subhead}%
%
%
\newcount\procno%
\global\procno=0%
\def\procinfo{\subheadinfo.\the\procno}%
\def\nextprocno{\global\advance\procno by 1 \procinfo}%
\makecounter{proc}%
%
%
\newcount\figno%
\global\figno=0%
\def\figinfo{\subheadinfo.\the\figno}%
\def\nextfigno{\global\advance\figno by 1 \figinfo}%
\makecounter{fig}%
%
%
\newcount\eqnno%
\global\eqnno=0%
\def\eqninfo{\text{{\rm (\headinfo.\the\eqnno)}}}%
\def\nexteqnno[#1]{\global\advance\eqnno by 1 \eqninfo\hbox{\eqnlabel{#1}}}%
\makecounter{eqn}%
%
%
\newcount\rmkno%
\global\rmkno=0%
\def\rmkinfo{\text{\subheadinfo.\the\rmkno}}%
\def\nextrmkno[#1]{\global\advance\rmkno by 1 \rmkinfo\hbox{\rmklabel{#1}}}%
\makecounter{rmk}%
%
%
%
%
%
\def\gobbleeight#1#2#3#4#5#6#7#8{}%
\newcount\citationno%
\global\citationno=0%
\def\citationinfo{\the\citationno}%
\makecounter{citation}%
\newwrite\biblio%
\def\newref#1#2{%
\def\temptext{#2}%
\edef\bibliotextoutput{\expandafter\gobbleeight\meaning\temptext}%
\global\advance\citationno by 1\citationlabel{#1}%
\ifmakebiblio%
    \edef\fileoutput{\write\biblio{\noindent\hbox to 0pt{\hss$[\the\citationno]$}\hskip 0.2em\bibliotextoutput\medskip}}%
    \fileoutput%
\fi}%
\def\cite#1{%
$[\citationref{#1}]$%
\ifmakebiblio%
    \edef\fileoutput{\write\biblio{#1}}%
    \fileoutput%
\fi%
}%
%
%
%
%
\let\mypar=\par%
\edef\Pagetitle={Blank}\headline={\hfil\Pagetitle\hfil}%
\edef\Pagefooter={Blank}\footline={\hfil\Pagefooter\hfil}%
%
%
\newcount\showpagenumflag%
\global\showpagenumflag=0 %
\def\nextoddpage%
{\newpage\ifodd\pageno%
\else\global\showpagenumflag=0 %
\null\vfil\eject%
\global\showpagenumflag=1 %
\fi}%
%
%
\font\headfont=cmb12%
\def\newhead#1[#2]%
{\ifhmode\mypar\fi%
\ifnum\headno=0 \else\goodbreak\bigskip\fi%
{\headfont\noindent\nextheadno\ - #1.}\headlabel{#2}%
\nobreak\medskip}%
%
%
\def\newsubhead#1[#2]%
{\ifhmode\mypar\fi%
\ifnum\subheadno=0 \else\goodbreak\medskip\fi%
{\bf\noindent\nextsubheadno\ - #1.\ }\subheadlabel{#2}}%
%
%
\newif\ifinproclaim%
\global\inproclaimfalse%
\def\proclaim#1{%
\goodbreak\medskip
\bgroup\inproclaimtrue%
\noindent{\bf #1}%
\nobreak\medskip\sl}%
\def\noskipproclaim#1{%
\goodbreak\medskip%
\bgroup\inproclaimtrue%
\noindent{\bf #1}\nobreak\sl}%
\def\endproclaim{\mypar\egroup\nobreak\medskip\ignorespaces}%
%
%
%
\newcount\xpos\newcount\ypos
\def\makelabelgrid{%
\xpos=-5 \ypos=-5 %
\loop\ifnum\xpos<6 %
{\loop\ifnum\ypos<6 %
\def\labeltext{x}%
\ifnum\xpos=0\def\labeltext{+}\fi%
\ifnum\ypos=0\def\labeltext{+}\fi%
\placelabel[\xpos][\ypos]{\labeltext}%
\advance\ypos by 1 %
\repeat}%
\advance\xpos by 1 %
\repeat}%
\def\placelabel[#1][#2]#3{{%
\setbox10=\hbox{\raise #2cm \hbox{\hskip #1cm #3}}%
\ht10=0pt \dp10=0pt \wd10=0pt \box10}}%
%
%
%
%
\def\myitem#1{\noindent\hbox to .5cm{\hfill#1\hss}}%
%
%
%
%
%
%
%
%
%
\font\sansseriften=cmss10%
\font\sansserifseven=cmss7%
\font\sansseriffive=cmss5%
\newfam\sansseriffam%
\textfont\sansseriffam=\sansseriften%
\scriptfont\sansseriffam=\sansserifseven%
\scriptscriptfont\sansseriffam=\sansseriffive%
\def\mathsf{\fam\sansseriffam}%
%
%
%
\font\boldten=cmb10%
\font\boldseven=cmb7%
\font\boldfive=cmb5%
\newfam\mathboldfam%
\textfont\mathboldfam=\boldten%
\scriptfont\mathboldfam=\boldseven%
\scriptscriptfont\mathboldfam=\boldfive%
\def\mathbf{\fam\mathboldfam}%
%
%
%
\font\mycmmiten=cmmi10%
\font\mycmmiseven=cmmi7%
\font\mycmmifive=cmmi5%
\newfam\mycmmifam%
\textfont\mycmmifam=\mycmmiten%
\scriptfont\mycmmifam=\mycmmiseven%
\scriptscriptfont\mycmmifam=\mycmmifive%
\def\hexa#1{\ifcase #1 0\or 1\or 2\or 3\or 4\or 5\or 6\or 7\or 8\or 9\or A\or B\or C\or D\or E\or F\fi}%
\mathchardef\mathi="7\hexa\mycmmifam7B%
\mathchardef\mathj="7\hexa\mycmmifam7C%
%
%
\font\mymsbmten=msbm10 at 8pt%
\font\mymsbmseven=msbm7 at 5.6pt
\font\mymsbmfive=msbm5 at 4pt%
\newfam\mymsbmfam%
\textfont\mymsbmfam=\mymsbmten%
\scriptfont\mymsbmfam=\mymsbmseven%
\scriptscriptfont\mymsbmfam=\mymsbmfive%
\mathchardef\mybeth="7\hexa\mymsbmfam69%
\mathchardef\mygimmel="7\hexa\mymsbmfam6A%
\mathchardef\mydaleth="7\hexa\mymsbmfam6B%
%
%
%
%
\def\proof{{\noindent\bf Proof:\ }}%
\def\remark[#1]{{\noindent\bf Remark \nextrmkno[#1].}}%
\def\qed{~$\square$}%
\def\makeop#1{\global\expandafter\def\csname op#1\endcsname{{\text{#1}}}}%
\def\makeopsmall#1{\global\expandafter\def\csname op#1\endcsname{{\text{\lowercase{#1}}}}}%
%
%
%
%
%
%
\makeop{Ext}%
\makeop{Int}%
\makeop{Dist}%
\makeop{Diam}%
\makeop{Length}%
%
%
%
%
%
\def\mlim{\mathop{{\text{Lim}}}}%
\def\mlimsup{\mathop{{\text{LimSup}}}}%
\def\mliminf{\mathop{{\text{LimInf}}}}%
\def\msup{\mathop{{\text{Sup}}}}%
\def\minf{\mathop{{\text{Inf}}}}%
%
%
\makeop{Dim}%
\makeop{Ker}%
\makeop{Coker}%
\makeop{Tr}%
\makeop{Adj}%
\makeop{Det}%
\makeop{End}%
\makeop{Lin}%
\makeop{Symm}%
\makeop{Mult}%
%
%
\makeop{dx}%
\makeop{dy}%
\makeop{dz}%
\makeop{dt}%
\makeop{dVol}%
\makeop{dArea}%
\makeop{Supp}%
\makeop{Hess}%
\makeop{Lip}%
%
%
\makeop{Re}%
\makeop{Im}%
\makeop{Arg}%
\makeop{Log}%
\makeop{Exp}%
%
%
\makeopsmall{Cos}%
\makeopsmall{Sin}%
\makeopsmall{Tan}%
\makeopsmall{Sec}%
\makeopsmall{Cosec}%
\makeopsmall{Cot}%
\makeopsmall{ArcCos}%
\makeopsmall{ArcSin}%
\makeopsmall{ArcTan}%
\makeopsmall{ArcSec}%
\makeopsmall{ArcCosec}%
\makeopsmall{ArcCot}%
%
%
\makeopsmall{Cosh}%
\makeopsmall{Sinh}%
\makeopsmall{Tanh}%
\makeopsmall{ArcCosh}%
\makeopsmall{ArcSinh}%
\makeopsmall{ArcTanh}%
%
%
\makeop{Vol}%
\makeop{Area}%
\makeop{Riem}%
\makeop{Ric}%
\makeop{Scal}%
\makeop{Euc}%
\makeop{Imm}%
\makeop{Emb}%
%
%
\makeop{Id}%
\makeop{Ad}%
\makeop{O}%
\makeop{SO}%
\makeop{SL}%
\makeop{GL}%
\makeop{Conf}%
\makeop{Homeo}%
\makeop{Diff}%
\makeop{Isom}%
%
%
\makeop{Ind}%
\makeop{Sig}%
\makeop{Spec}%
%
%
\makeop{Conv}%
\makeop{Max}%
\makeop{Min}%
\makeop{Mod}%
\makeop{Deg}%
\makeop{loc}%
%
%
%
%
%
%
%
%
%
%
%
%
%
 %
%
%
%
%
%
\def\Pagetitle{\hfil\ifnum\pageno=1\null\else{\rm On the Weyl problem in Minkowski space.}\fi\hfil}
\def\Pagefooter{\hfil{\myfontdefault\folio}\hfil}
\font\tablefont=cmr7
\newif\ifshowaddress\showaddresstrue
\def\centre{\rightskip=0pt plus 1fil \leftskip=0pt plus 1fil \spaceskip=.3333em \xspaceskip=.5em \parfillskip=0em \parindent=0em}%
\def\textmonth#1{\ifcase#1\or January\or Febuary\or March\or April\or May\or June\or July\or August\or September\or October\or November\or December\fi}
\font\abstracttitlefont=cmr10 at 12pt {\abstracttitlefont\centre On the Weyl problem in Minkowski space.\par}
\bigskip
{\centre 8th June 2018\par}
\bigskip
{\centre Graham Smith\footnote{${}^*$}{{\tablefont Instituto de Matem\'atica, UFRJ, Av. Athos da Silveira Ramos 149, Centro de Tecnologia - Bloco C, Cidade Universit\'aria - Ilha do Fund\~ao, Caixa Postal 68530, 21941-909, Rio de Janeiro, RJ - BRAZIL\hfill}}\par}
\bigskip
\noindent{\bf Abstract:~}Let $S$ be a closed surface of hyperbolic type. We show that, for every pair $(g_+,g_-)$ of negatively curved metrics over $S$ there exists a unique GHMC Minkowski spacetime $X$ into which $(S,g_+)$ and $(S,g_-)$ isometrically embed as Cauchy surfaces in the future and past components respectively.
%
%
\bigskip
\myfontdefault
\catcode`\@=11
\def\triplealign#1{\null\,\vcenter{\openup1\jot \m@th %
\ialign{\strut\hfil$\displaystyle{##}\quad$&$\displaystyle{{}##}$\hfil&$\displaystyle{{}##}$\hfil\crcr#1\crcr}}\,}
\def\multiline#1{\null\,\vcenter{\openup1\jot \m@th %
\ialign{\strut$\displaystyle{##}$\hfil&$\displaystyle{{}##}$\hfil\crcr#1\crcr}}\,}
\catcode`\@=12

\makeop{ad}
\makeop{PSO}
\makeop{PSL}
\makeop{II}
\makeop{III}
\makeop{Width}
\makeop{hyp}
\makeop{rep}
\makeop{T}
\makeop{SGr}
\makeop{Coth}
\makeop{GHMC}
\makeop{EC}
\makeop{Lam}
\makeop{bdd}
\makeop{A}
\makeop{D}
\makeop{E}
\makeop{F}
\makeop{G}
\makeop{J}
\makeop{R}
\makeop{Teich}
\makeop{L}
\makeop{dArea}

\makeop{Mod}
\def\opConf{{\text{Conf}}}
\makeop{Inj}
\def\optildeConf{{\text{C}\widetilde{\text{onf}}}}
\def\optildeHyp{{\text{H}\widetilde{\text{yp}}}}
\makeop{Hyp}
\makeop{Diff}
\makeop{h}
\makeop{Sys}
\makeop{T}
\def\Surface{S}
\makeop{QF}
\makeop{AD}
\newref{Ahlfors}{Ahlfors L., {\sl Lectures on quasiconformal mappings}, AMS, (2006)}
\newref{AndEtAl}{Andersson L., Barbot T., Benedetti R., Bonsante F., Goldman W. M., Labourie F., Scannell K. P., Schlenker J. M., Notes on: “Lorentz spacetimes of constant curvature”, {\sl Geom. Dedicata}, {\bf 126}, (2007), 47--70}
\newref{BarbFill}{Barbot T., Fillastre F., Quasi-Fuchsian co-Minkowski manifolds, to appear {\sl In the tradition of Thurston}, Springer-Verlag, (2020)}
\newref{Bers}{Bers L., Simultaneous uniformization, {\sl Bull. Amer. Math. Soc.}, {\bf 66}, no.2, (1960), 94--97}
\newref{BonMonSch}{Bonsante F., Mondello G., Schlenker J. M., A cyclic extension of the earthquake flow II, {\sl Ann. Sci. Ec. Norm. Sup\'er.}, {\bf 48}, (2015), no. 4, 811--859}
\newref{BonSeppi}{Bonsante F., Seppi A., On Codazzi tensors on a hyperbolic surface and flat Lorentzian geometry, {\sl Int. Math. Res. Not. IMRN}, (2):343–417, (2016)}
\newref{Buser}{Buser P., {\sl Geometry and spectra of compact Riemann surfaces}, Birkh\"auser Verlag, (1992)}
\newref{FillSmi}{Fillastre F., Smith G., Group actions and scattering problems in Teichm\"uller theory, in {\sl Handbook of group actions IV}, Advanced Lectures in Mathematics, {\bf 40}, (2018), 359--417}
\newref{FillSmiII}{Fillastre F., Smith G., A note on invariant constant curvature immersions in Minkowski space, to appear in {\sl Geom. Dedicata}}
\newref{GilbTrud}{Gilbarg D., Trudinger N. S., {\sl Elliptic partial differential equations of second order}, Classics in Mathematics, Springer-Verlag, Berlin, (2001)}
\newref{Jost}{Jost J., {\sl Riemannian Geometry and Geometric Analysis}, Universitext, Springer-Verlag, (2011)}
\newref{LabI}{Labourie F., Probl\`eme de Minkowski et surfaces \`a courbure constante dans les vari\'et\'es hyperboliques, {\sl Bull. Soc. math. France}, {\bf 119}, (1991), 307--325}
\newref{LabII}{Labourie, F., Metriques prescrites sur le bord des vari\'et\'es hyperboliques de dimension $3$, {\sl J. Differ. Geom.}, {\bf 35}, (1992), 609--626}\newref{Mess}{Mess G., Lorentz spacetimes of constant curvature, {\sl Geom. Dedicata}, {\bf 126}, (2007), 3--45}
\newref{Mumford}{Mumford D., A remark on Mahler's compactness theorem, {\sl Proc. AMS}, {\bf 28}, no.1, (1971), 289--294}
\newref{NewNir}{Newlander A., Nirenberg L., Complex analytic coordinates in almost complex manifolds, {\sl Ann. of Math.}, {\bf 65}, (1957), 391--404}
\newref{Schlenker}{Schlenker J. M., Hyperbolic manifolds with convex boundary, {\sl Inventiones mathematicae}, {\bf 163}, (2006), 109--169}
\newref{Schoen}{Schoen R., The role of harmonic mappings in rigidity and deformation problems, in {\sl Complex geometry: Proceedings of the Osaka international conference}, (1993), 179-200}
\newref{SchoenYau}{Schoen R., Yau S. T., On univalent harmonic maps between surfaces, {\sl Invent. Math.}, {\bf 44}, (1978), 265--278}
\newref{Tamburelli}{Tamburelli A., Prescribing metrics on the boundary of anti-de Sitter $3$-manifolds, {\sl Int. Math. Res. Not. IMRN}, (2018), no. 5, 1281--1313}
\newref{Thurston}{Thurston W., The geometry and topology of three manifolds}
\newref{TrapaniValli}{Trapani S., Valli G., One-harmonic maps on Riemann surfaces, {\sl Comm. Anal. Geom.}, {\bf 3}, no. 4, (1985), 645--681}
\newref{Tromba}{Tromba A. J., {\sl Teichm\"uller theory in Riemannian geometry}, Lectures in Mathematics, ETH Z\"urich, Birkh\"auser Verlag, Basel, (1992)}
\newref{Wolf}{Wolf M., The Teichm\"uller theory of harmonic maps, {\sl J. Diff. Geom}, {\bf 29}, (1989), 449--479}
\newhead{Introduction}[Introduction]
\newsubhead{Equivariant isometric embeddings}[EquivariantIsometricEmbeddings]
The study of spaces of representations of surface groups in certain types of Lie groups has become a central theme of modern Teichm\"uller theory. In this paper, we will be concerned with the prescription of such representations by geometric data. Consider a closed surface $S$ of hyperbolic type and denote its fundamental group by $\Gamma$. The simultaneous uniformisation theorem proven by Bers in \cite{Bers} can be viewed as a smooth parametrisation of the space $\opQF(\Gamma,\opPSO(3,1))$ of conjugacy classes of quasi-Fuchsian representations of $\Gamma$ in $\opPSO(3,1)$ by pairs of marked conformal structures over $S$. In a similar vein, the work \cite{LabI} of Labourie yields a two-dimensional family of smooth parametrisations of the same space by pairs of marked hyperbolic metrics over $S$. Furthermore, Labourie's family interpolates between Bers' parametrisation and another parametrisation, also by pairs of marked hyperbolic metrics, conjectured by Thurston in \cite{Thurston}. In the case where the target group is $\opPSO(2,2)$, a one-parameter family of smooth parametrisations of the space $\opQF(\Gamma,\opPSO(2,2))$ of conjugacy classes of quasi-Fuchsian representations of $\Gamma$ is constructed in the work \cite{BonMonSch} of Bonsante, Mondello \& Schlenker. Likewise, in the case where the target group is $\opPSO(2,1)\ltimes\Bbb{R}^{2,1}$, a two-parameter family of smooth parametrisations of the space $\opAD(\Gamma,\opPSO(2,1)\ltimes\Bbb{R}^{2,1})$ of conjugacy classes of affine deformations (see below) is constructed in our joint work \cite{FillSmiII} with Fran\c{c}ois Fillastre. These and similar parametrisations are discussed in greater detail in \cite{FillSmi}.
\par
A related problem is that of prescribing certain types of representations by more general geometric data. For example, by proving existence of unique solutions to a certain non-linear PDE, the respective results \cite{LabII} and \cite{Schlenker} of Labourie and Schlenker together show how elements of $\opQF(\Gamma,\opPSO(3,1))$ are determined uniquely by pairs of negatively curved metrics over $S$. In a similar vein, in \cite{Tamburelli}, by proving existence of solutions to a different non-linear PDE, Tamburelli associates elements of $\opQF(\Gamma,\opPSO(2,2))$ to arbitrary pairs of metrics over $S$ of curvature bounded above by $-1$. However, uniqueness of the representations constructed by Tamburelli remains an open problem.
\par
In the current paper, extending our joint work \cite{FillSmiII} with Fran\c{c}ois Fillastre, we will be concerned with the prescription of elements of $\opAD(\Gamma,\opPSO(2,1)\ltimes\Bbb{R}^{2,1})$ by pairs of negatively curved metrics over $S$. We first recall some algebraic preliminaries. Let $\Bbb{R}^{2,1}$ be $(2+1)$-dimensional Minkowski space with metric given by
$$
\delta^{2,1}(x,y) := x_1y_1 + x_2y_2 - x_3y_3.\eqnum{\nexteqnno[MinkowskiMetric]}
$$
Let $\opO(2,1)$ be its group of linear isometries. Its group of affine isometries is given by the semi-direct product $\opO(2,1)\ltimes\Bbb{R}^{2,1}$ with group law
$$
(g,x)\cdot(h,y) := (gh,x + gy).\eqnum{\nexteqnno[SemiDirectGroupLaw]}
$$
Let $\opSO_0(2,1)$ denote the identity component of $\opO(2,1)$. This subgroup consists of precisely those isometries which preserve both the orientation and the time orientation. Trivially $\opSO_0(2,1)\ltimes\Bbb{R}^{2,1}$ is also the identity component of $\opO(2,1)\ltimes\Bbb{R}^{2,1}$.
\par
A homomorphism $\rho:\Gamma\rightarrow\opSO_0(2,1)$ is said to be {\sl Fuchsian} whenever it is injective with discrete image. Given a homomorphism $(\rho,\tau):\pi_1(S)\rightarrow\opSO_0(2,1)\ltimes\Bbb{R}^{2,1}$, the functions $\rho$ and $\tau$ are called its {\sl linear} and {\sl cocycle} components respectively. The homomorphism $(\rho,\tau)$ is said to be an {\sl affine deformation} whenever its linear component $\rho$ is Fuchsian. Recall that the space of conjugacy classes of Fuchsian homomorphisms identifies with the Teichm\"uller space $\opTeich[S]$ of marked hyperbolic metrics over $S$. Likewise, the space $\opAD(\Gamma,\opSO(2,1)\ltimes\Bbb{R}^{2,1})$ of conjugacy classes of affine deformations identifies with the total space of the tangent bundle $\opT\opTeich[S]$ of Teichm\"uller space (see \cite{BarbFill} and \cite{BonSeppi}).
\par
We also require some terminology concerning immersions in Minkowski space. Let $\tilde{S}$ be the universal cover of $S$. A smooth immersion $e:\tilde{S}\rightarrow\Bbb{R}^{2,1}$ is said to be {\sl spacelike} whenever its first fundamental form is positive-definite. Every spacelike immersion has a well-defined, future-oriented, unit normal vector field $N$. A spacelike immersion $e$ is said to {\sl locally strictly convex} (LSC) whenever its second fundamental form $\opII$ with respect to this normal is either positive- or negative-definite. The immersion is itself is said to be {\sl future-oriented} whenever $\opII$ is positive-definite and {\sl past-oriented} otherwise. Finally, given an affine deformation $(\rho,\tau):\Gamma\rightarrow\opSO_0(2,1)\ltimes\Bbb{R}^{2,1}$, the immersion $e$ is said to be $(\rho,\tau)$-equivariant whenever it satisfies, for all $x\in\tilde{S}$ and for all $\gamma\in\Gamma$,
$$
e(\gamma x) = \rho(\gamma)e(x) + \tau(\gamma).\eqnum{\nexteqnno[Equivariance]}
$$
We prove
\proclaim{Theorem \nextprocno}
\noindent Let $g_+$ and $g_-$ be negatively-curved metrics over $S$. There exists an affine deformation $(\rho,\tau):\Gamma\rightarrow\opSO_0(2,1)\ltimes\Bbb{R}^{2,1}$, a future-oriented, LSC, $(\rho,\tau)$-equivariant spacelike immersion $e_+:\tilde{S}\rightarrow\Bbb{R}^{2,1}$ and a past-oriented, LSC, $(\rho,\tau)$-equivariant spacelike immersion $e_-:\tilde{S}\rightarrow\Bbb{R}^{2,1}$ such that
$$
g_\pm = e_\pm^*\delta^{2,1}.
$$
Furthermore, $(\rho,\tau)$ is unique up to conjugation by an element of $\opSO_0(2,1)\ltimes\Bbb{R}^{2,1}$ and, given $(\rho,\tau)$, $e_+$ and $e_-$ are also unique.
\endproclaim
\proclabel{MainTheoremI}
\remark[MainTheoremI] Theorem \procref{MainTheoremI} is proven in Section \subheadref{ProofOfMainTheorem}, below.
\newsubhead{GHMC Minkowski spacetimes}[GHMCMinkowskiSpacetimes]
We now reformulate Theorem \procref{MainTheoremI} in terms of the prescription of certain types of semi-riemannian manifolds by pairs of metrics over $S$. We first recall some terminology of semi-riemannian geometry (c.f. \cite{FillSmi}). A {\sl Minkowski spacetime} is a semi-riemannian manifold which is everywhere locally isometric to $\Bbb{R}^{2,1}$. A smoothly embedded curve in a Minkowski spacetime is said to be {\sl causal} whenever its derivative has non-positive norm-squared at every point. A Minkowski spacetime is itself said to be {\sl causal} whenever it contains no non-trivial, closed, causal curve. A causal spacetime is said to be {\sl globally hyperbolic} whenever it contains a {\sl Cauchy hypersurface}, that is, a smoothly embedded hypersurface that meets every inextensible, causal curve exactly once. A globally hyperbolic spacetime $X$ is said to be {\sl maximal} whenever it cannot be isometrically embedded into a strictly larger globally hyperbolic spacetime $X'$ in such a manner that the Cauchy hypersurfaces of $X$ are mapped to Cauchy hypersurfaces of $X'$. Finally, a globally hyperbolic spacetime is said to be {\sl Cauchy compact} whenever its Cauchy hypersurface, which is unique up to diffeomorphism, is compact. A Minkowski spacetime which possesses all the above properties is said to be {\sl GHMC} (Globally, Hyperbolic, Maximal and Cauchy Compact).
\par
In \cite{Mess} (c.f. \cite{AndEtAl} and \cite{FillSmi}), Mess proves a key structure theorem relating GHMC Minkowski spacetimes to affine deformations which we now recall. Let $(\rho,\tau):\Gamma\rightarrow\opSO_0(2,1)\ltimes\Bbb{R}^{2,1}$ be an affine deformation. There exists a unique closed, convex, proper subset $K_+\subset\Bbb{R}^{2,1}$ such that
\medskip
\myitem{(1)} $K_+$ is future-complete in the sense that every future-oriented causal curve starting in $K_+$ remains in $K+$,
\medskip
\myitem{(2)} $K_+$ is invariant under the action of $(\rho,\tau)$, and
\medskip
\myitem{(3)} $K_+$ is maximal amongst all such subsets in the sense that if $K_+'$ is another closed, convex, proper subset of $\Bbb{R}^{2,1}$ satisfying the above two conditions, then $K_+'\subseteq K_+$.
\medskip
\noindent There likewise exists a unique closed, convex, proper subset $K_-\subseteq\Bbb{R}^{2,1}$ which is past-complete in the obvious sense, invariant and maximal. The sets $K_+$ and $K_-$ have disjoint interiors and $(\rho,\tau)$ acts properly discontinuously over these interiors. The quotients $K_+/(\rho,\tau)$ and $K_-/(\rho,\tau)$ are respectively the future- and past-complete components of a marked GHMC Minkowski spacetime with Cauchy hypersurface diffeomorphic to $S$. In fact, all such GHMC Minkowski spacetimes arise in this manner, and two affine deformations yield the same GHMC Minkowski spacetime if and only if they are conjugate by an element of $\opSO_0(2,1)\ltimes\Bbb{R}^{2,1}$. In this manner, Mess obtains an identification of the space $\opGHMC_0[S]$ of marked GHMC Minkowski spacetimes with Cauchy hypersurface diffeomorphic to $S$ with the space $\opT\opTeich[S]$ of affine deformations of $\Gamma$.
\par
Theorem \procref{MainTheoremI} is now expressed in terms of the Weyl problem concerning the construction of isometric embeddings in GHMC Minkowski spacetimes.
\proclaim{Theorem \nextprocno}
\noindent Given two negatively curved metrics $g_\pm$ over $S$, there exists a unique GHMC Minkowski spacetime into which $(S,g_+)$ and $(S,g_-)$ embed isometrically as Cauchy surfaces in its future and past components respectively. Furthermore, these embeddings are also unique.
\endproclaim
\newsubhead{One-harmonic maps}[OneHarmonicMaps]
Our proof rests on the work \cite{TrapaniValli} of Trapani \& Valli which we now briefly review. Let $\opDiff(\Surface)$ denote the group of smooth diffeomorphisms of $S$. Given two metrics $g$ and $h$ over $\Surface$, Trapani \& Valli define the real-valued functional $T_\partial(g,h,\cdot)$ over $\opDiff(\Surface)$ by
$$
T_\partial(g,h,\Phi) := \int_\Surface\|\partial\Phi\|\opdArea_g,\eqnum{\nexteqnno[TVEnergyFunctional]}
$$
where $\opdArea_g$ here denotes the area form of $g$, $\partial\Phi$ denotes the complex linear component of $D\Phi$ and $\|\partial\Phi\|$ denotes its operator norm with respect to $g$ and $h$. In what follows, we will call this functional the {\sl $(1,0)$-energy} of $\Phi$. Following Trapani \& Valli, its critical points are called {\sl one-harmonic} diffeomorphisms. Trapani \& Valli show
\proclaim{Theorem \nextprocno, {\bf Trapani \& Valli (1995)}}
\noindent For any pair $(g,h)$ of marked, negatively-curved metrics over $\Surface$, there exists a unique one-harmonic diffeomorphism $\Phi:(S,g)\rightarrow (S,h)$ which preserves the marking. Furthermore, $\Phi$ depends smoothly on $g$ and $h$.
\endproclaim
\proclabel{TrapaniValliI}
In order to apply Theorem \procref{TrapaniValliI} in the present context, we first reformulate it in the language of Codazzi tensors which we now recall (c.f. \cite{FillSmi} and \cite{Tromba}). Let $g$ be a smooth metric over $S$ and let $\nabla$ be its Levi-Civita covariant derivative. Let $\opEnd(T\Surface)$ be the bundle of endomorphisms of $T\Surface$. Sections of this bundle are called {\sl endomorphism fields}. An endomorphism field $A$ is said to be a {\sl Codazzi field} of $g$ whenever it is symmetric with respect to $g$ and
$$
d^\nabla A = 0.\eqnum{\nexteqnno[CodazziFieldCondition]}
$$
In Corollary \procref{CorDerivativeOfEnergyAtSym} below, we show that a smooth diffeomorphism $\Phi:(S,g)\rightarrow (S,h)$ is one-harmonic if and only if there exists a positive-definite Codazzi field $A$ of $g$ such that
$$
\Phi^* h = g(A\cdot,A\cdot).
$$
Theorem \procref{TrapaniValliI} therefore becomes
\proclaim{Theorem \nextprocno, {\bf Trapani \& Valli (1995)}}
\noindent For any pair $(g,h)$ of marked, negatively-curved metrics over $\Surface$, there exists a unique pair $(A,\Phi)$, where $A$ is a positive-definite Codazzi field of $g$ and $\Phi:(S,g)\rightarrow(S,h)$ is a smooth diffeomorphism which preserves the marking, such that
$$
\Phi^*h = g(A\cdot,A\cdot).
$$
Furthermore, both $A$ and $\Phi$ vary smoothly with $g$ and $h$.
\endproclaim
\proclabel{TrapaniValliII}
Theorem \procref{TrapaniValliI} is proven in \cite{TrapaniValli} using the continuity method which, we recall, consists of two parts, namely, a local perturbation result and a compactness result. In studying Trapani \& Valli's work, we have found that the formulae required to understand the perturbation part become much simpler when expressed in the language of Codazzi fields. Indeed, this is particularly true for the generalized lorentzian metric and its geodesics studied in Section $5$ of \cite{TrapaniValli} (see Appendix \headref{TheGeneralisedLorentzianMetric}, below). With the aim of bringing Trapani \& Valli's work to a wider audience, Sections \headref{Preliminaries} and \headref{PerturbationTheory} of this paper are devoted to reformulating the proof of this part of their result in the language of Codazzi fields. Since the details of this argument are already presented in \cite{TrapaniValli}, we will only present in these sections the main ideas, providing only those proofs we consider necessary.
\par
In Section \headref{Compactness}, we modify the proof of the compactness part of \cite{TrapaniValli}. In particular, by clarifying the role played by each of the two metrics, we obtain the following new compactness result for sequences of positive-definite Codazzi fields, which we will require in Section \headref{TheWeylProblem} and which we believe to be of independent interest.
\proclaim{Theorem \nextprocno}
\noindent Let $(g_m)_{m\in\Bbb{N}}$ and $(\phi_m)_{m\in\Bbb{N}}$ be respectively sequences of smooth metrics and smooth, positive functions over $\Surface$ converging respectively in the $C^\infty$ sense to the smooth metric $g_\infty$ and the smooth, positive function $f_\infty$. For all $m$, let $A_m$ be a positive-definite Codazzi field of $g_m$ such that $\opDet(A_m)=\phi_m$. If
$$
\msup_m\int_\Surface\opTr(A_m)\opdArea_{g_m} < \infty,
$$
then there exists a positive-definite Codazzi field $A_\infty$ of $g_\infty$ towards which $(A_m)_{m\in\Bbb{N}}$ subconverges in the $C^\infty$ sense.
\endproclaim
\proclabel{IntroCompactness}
\remark[IntroCompactness] Theorem \procref{IntroCompactness} is proven in Theorem \procref{CompactnessForCodazziFields}, below.
\medskip
\remark[IntroCompactnessII] The proof of compactness presented in \cite{TrapaniValli}, although correct, is confused. In the notation of that paper, it is not necessary to prove that the conformal classes of the metrics $(m_n)$ subconverge (Lemma $6.3$). The subsequent results concerning the properties of this sequence are also unnecessary. However, it is necessary to prove that the conformal classes of the metrics $(p_n)$ subconverge. Although this readily follows from Lemmas $6.3$ and $6.6$, this is not mentioned explicitely.
\newsubhead{Acknowledgements}[Acknowledgements]
The author is grateful to Fran\c{c}ois Fillastre for invaluable contributions throughout the preparation of this paper. We are also grateful to Francesco Bonsante for having drawn our attention to the work \cite{TrapaniValli} of Trapani \& Valli.
\newhead{Preliminaries}[Preliminaries]
\newsubhead{Linear algebra}[LinearAlgebraicPreliminaries]
In this and the following chapter, we reformulate the perturbation part of \cite{TrapaniValli} using the terminology of Codazzi fields. We underline that, since the results presented in these chapters have already been proven in that paper, we provide only those proofs that we consider necessary.
\par
Let $\opEnd(2)$ denote the space of real $2\times 2$ matrices. Let $\opSymm_+(2)$ denote the space of positive-definite, symmetric bilinear forms over $\Bbb{R}^2$. Let $g\in\opSymm_+(2)$ denote the standard metric. We will study $\opSymm_+(2)$ as a symmetric space. Let $\opGL_+(2)$ denote the group of orientation-preserving $2\times 2$ matrices and let $\opSO(2)$ denote its special-orthogonal subgroup. The action of $\opGL_+(2)$ on $\opSymm_+(2)$ is given by
$$
(Ah)(\cdot,\cdot) := h(A\cdot,A\cdot).\eqnum{\nexteqnno[ActionOfGL]}
$$
This action is transitive and its stabiliser of $g$ is $\opSO(2)$ so that $\opSymm_+(2)$ identifies with the symmetric space $\opSO(2)\backslash\opGL_+(2)$. In particular, $\opGL_+(2)$ is a principal $\opSO(2)$-bundle over $\opSymm_+(2)$ with projection $h:\opGL_+(2)\rightarrow\opSymm_+(2)$ given by
$$
h(A) := g(A\cdot,A\cdot).\eqnum{\nexteqnno[BundleProjection]}
$$
Conversely, let $A:\opSymm_+(2)\rightarrow\opGL_+(2)$ be such that, for all $h$, $A:=A(h)$ is the unique positive-definite, symmetric matrix such that $h=h(A)$. This function defines a global section of $\opGL_+(2)$ over $\opSymm_+(2)$.
\par
We now introduce the $(1,0)$-seminorm over $\opEnd(2)$ which is one of the main concepts used in \cite{TrapaniValli}. The standard inner-product over $\opEnd(2)$ is defined by
$$
\langle A,B\rangle := \opTr(AB^t).\eqnum{\nexteqnno[MatrixInnerProduct]}
$$
Let
$$
J := \pmatrix 0\hfill& -1\hfill\cr 1\hfill& 0\hfill\cr\endpmatrix\eqnum{\nexteqnno[DefnComplexStructure]}
$$
denote the standard complex structure of $\Bbb{R}^2$. Given a matrix $A\in\opEnd(2)$, its $(1,0)$-{\sl seminorm} is defined by
$$
\sigma_g(A) := \|A^{1,0}\|,\eqnum{\nexteqnno[DefnOneZeroSeminorm]}
$$
where $\|\cdot\|$ here denotes the norm of the standard inner-product \eqnref{MatrixInnerProduct} and
$$
A^{1,0} := \frac{1}{2}(A -  J A J)\eqnum{\nexteqnno[JLinAndAntilin]}
$$
is the $J$-linear component of $A$. It is straightforward to show that
$$
A^{1,0} = \frac{1}{2}\opTr(A)\opId - \frac{1}{2}\opTr(AJ)J,\eqnum{\nexteqnno[FormulaForComplexLinearComponent]}
$$
from which we obtain the useful formula
$$
\sigma_g(A)^2 = \frac{1}{2}\opTr(A)^2 + \frac{1}{2}\opTr(JA)^2.\eqnum{\nexteqnno[FormulaForSigma]}
$$
In particular, for all $M\in\opSO(2)$ and for all $A\in\opEnd(2)$,
$$
\sigma_g(MA) = \sigma_g(A),\eqnum{\nexteqnno[InvarianceOfSigma]}
$$
so that $\sigma_g$ descends to a function defined over $\opSymm_+(2)$. This function will also be called the $(1,0)$-{\sl seminorm} and will also be denoted by $\sigma_g$. In the notation of \cite{TrapaniValli}, $\sigma_g(A)$ is equal to $\|\partial\opId\|$ with respect to the metrics $g$ and $h(A)$.
\par
Finally, observe that, for all $A\in\opEnd(2)$,
$$\eqalign{
A-A^t &= -\opTr(AJ)J\ \text{and}\cr
JA^tJ &= -\opDet(A)A^{-1}.}\eqnum{\nexteqnno[TwoDimMatrixRelations]}
$$
In particular, $A$ is symmetric if and only if $\opTr(AJ)$ vanishes. Differentiating \eqnref{FormulaForSigma} therefore yields
\proclaim{Lemma \nextprocno}
\noindent If $\opTr(A)\neq 0$, then the first derivative of $\sigma_g$ at $A$ satisfies
$$
D\sigma_g(A)\cdot B = \opTr(B) + \frac{\opTr( J A)}{\opTr(A)}\opTr( J B) + (F(A)\cdot B)\opTr(JA)^2,\eqnum{\nexteqnno[FirstDerivativeA]}
$$
for some smooth function $F$. In particular, when $A$ is symmetric and positive definite,
$$
D\sigma_g(A)\cdot B = \opTr(B).\eqnum{\nexteqnno[FirstDerivativeB]}
$$
\endproclaim
\proclabel{FirstDerivative}
\newsubhead{Calculus}[GradientAndDivergence]
Let $\Omega$ be an open subset of $\Bbb{R}^2$. Let $J$ be as in \eqnref{DefnComplexStructure}. Let $g:=\langle\cdot,\cdot\rangle$ now be a riemannian metric over $\Omega$ conformal with respect to $J$. Let $\nabla$ denote its Levi-Civita covariant derivative. We recall three different types of derivative constructed out of $\nabla$ that will be used in the sequel. First, let $e_1,e_2$ be an oriented, orthonormal frame of $g$. The {\sl gradient} of any function $f:\Omega\rightarrow\Bbb{R}$ is the vector field defined by
$$
\nabla f := \sum_{i=1}^2(D_{e_i}f)e_i,\eqnum{\nexteqnno[DefnGradientFunction]}
$$
the {\sl divergence} of any vector field $X:\Omega\rightarrow\Bbb{R}^2$ is the function defined by
$$
\nabla\cdot X := \sum_{i=1}^2\langle\nabla_{e_i}X,e_i\rangle\eqnum{\nexteqnno[DefnDivergenceVectorField]}
$$
and the {\sl divergence} of any endomorphism field $A:\Omega\rightarrow\opEnd(2)$ is the vector field defined by
$$
\nabla\cdot A := \sum_{i=1}^2(\nabla_{e_i}A)e_i.\eqnum{\nexteqnno[DefnDivergenceEndmorphism]}
$$
For all $f$, $X$ and $A$, the vector field $\nabla f$, the function $\nabla\cdot X$ and the vector field $\nabla\cdot A$ are independent of the orthonormal frame chosen. Furthermore, they are related to the exterior derivatives by
$$\eqalign{
\nabla f &= df^\#,\cr
\nabla\cdot X &= -d(\alpha J)(e_1,e_2)\ \text{and}\cr
\nabla\cdot A &= -d^\nabla(A J)(e_1,e_2),\cr}\eqnum{\nexteqnno[ObjectsInTermsOfExtAlg]}
$$
where here $\alpha:=X^\flat$ and $\#$ and $\flat$ denote Berger's musical isomorphisms.
\proclaim{Lemma \nextprocno}
\noindent For every endomorphism field $A$,
$$
\nabla\opTr(A) - \nabla\cdot A -  J\nabla\opTr(A J) +  J\nabla\cdot(A J) = 0.\eqnum{\nexteqnno[FundamentalRelation]}
$$
\endproclaim
\proclabel{FundamentalRelation}
\remark[FundamentalRelation] Before proving Lemma \procref{FundamentalRelation}, recall that, by \eqnref{TwoDimMatrixRelations}, an endomorphism $A$ is symmetric if and only if $AJ$ is trace-free. Thus, if an endomorphism field $A$ is trace-free, symmetric and closed in the sense that $\nabla\cdot(A J)=0$, it follows by \eqnref{FundamentalRelation} that it is also co-closed, and vice-versa. Classical Hodge theory then shows that the space of closed (or co-closed) trace-free, symmetric endomorphism fields over a given compact Riemann surface $\Sigma$ is finite-dimensional. Furthermore, when $\Sigma$ is of hyperbolic type, its dimension is equal to $(6g-6)$, where $g$ here denotes the genus of $\Sigma$. Observe finally that two natural complex structures are defined over this space by precomposition and postcomposition with $J$ and that these complex structures differ from one-another only by a change of sign.
\medskip
\proof Choose a point $x\in\Omega$. We may assume that the frame $(e_1,e_2)$ is parallel at $x$. Then, using \eqnref{TwoDimMatrixRelations} and \eqnref{ObjectsInTermsOfExtAlg} we obtain, for every vector field $X$,
$$\eqalign{
\langle\nabla\opTr(A),X\rangle
&=\sum_i\langle(\nabla_{X} A)e_i,e_i\rangle\cr
&=\sum_i\langle(\nabla_{e_i} A)X,e_i\rangle + \sum_i\langle(d^\nabla A)(X,e_i),e_i\rangle\cr
&=\sum_i\langle(\nabla_{e_i}A)e_i,X\rangle - \sum_i\langle(\nabla_{e_i}\opTr(AJ)J)X,e_i\rangle +\sum_i\langle (d^\nabla A)(X,e_i),e_i\rangle\cr
&=\langle\nabla\cdot A,X\rangle - \langle\nabla\opTr(AJ),JX\rangle + \langle\nabla\cdot(AJ),JX\rangle,\cr}
$$
and the result follows.\qed
\newhead{Perturbation theory}[PerturbationTheory]
\newsubhead{The $(1,0)$-energy density and the $(1,0)$-energy}[TheEnergyDensity]
We continue our reformulation of the perturbation part of \cite{TrapaniValli} using the terminology of Codazzi fields. We underline once again that we will only provide in this chapter those proofs that we consider necessary.
\par
Let $\Omega$ be an open subset of $\Bbb{R}^2$. Let $J$ be the standard complex structure of $\Bbb{R}^2$ as given in \eqnref{DefnComplexStructure}. Let $g:=\langle\cdot,\cdot\rangle$ be a riemannian metric over $\Omega$ conformal with respect to $J$. We extend the structures described in Section \subheadref{LinearAlgebraicPreliminaries} to smooth families fibred over $\Omega$. First, given any finite-dimensional vector space $E$, let $C^\infty_\opbdd(\Omega,E)$ denote the Frechet space of smooth functions $A:\Omega\rightarrow E$ all of whose derivatives to all orders are bounded. Let $\opSymm_+(2,\Omega)$ denote the open subset of $C^\infty_\opbdd(\Omega,\opSymm(\Bbb{R}^2))$ consisting of those functions $h$ having the property that there exists $C>0$ such that, for all $z\in\Omega$, and for all $\xi\in\Bbb{R}^2$,
$$
h(z)(\xi,\xi) \geqslant \frac{1}{C}\|\xi\|^2.
$$
We interpret $\opSymm_+(2,\Omega)$ as the space of suitably regular riemannian metrics over $\Omega$.
\par
As before, we study $\opSymm_+(2,\Omega)$ as a symmetric space. Let $\opGL_+(2,\Omega)$ denote the multiplicative group of functions $A\in C^\infty_\opbdd(\Omega,\opEnd(2))$ with the property that there exists $C>0$ such that, for all $z\in\Omega$,
$$
\opDet(A(z)) \geq \frac{1}{C}.
$$
Likewise, let $\opSO(2,\Omega)$ denote the subgroup consisting of those functions $M$ such that, for all $z\in\Omega$, $M(z)$ is special-orthogonal with respect to the metric $g(z)$. The action of $\opGL_+(2,\Omega)$ over $\opSymm_+(2,\Omega)$ is given by
$$
Ah := h(A\cdot,A\cdot).
$$
This action is transitive and its stabiliser of $g$ is $\opSO(2,\Omega)$ so that $\opSymm_+(2,\Omega)$ identifies with the symmetric space $\opSO(2,\Omega)\backslash\opGL_+(2,\Omega)$. In particular, $\opGL_+(2,\Omega)$ is a principal $\opSO(2,\Omega)$-bundle over $\opSymm_+(2,\Omega)$ with projection $\oph:\opGL_+(2,\Omega)\rightarrow\opSymm_+(2,\Omega)$ given by
$$
\oph[A] := g(A\cdot,A\cdot).
$$
Conversely, define $\opA:\opSymm_+(2,\Omega)\rightarrow\opGL_+(2,\Omega)$ such that, for all $h$, $A:=\opA[h]$ is the unique positive-definite, symmetric matrix field such that $h=\oph[A]$. This function defines a global section of $\opGL_+(2,\Omega)$ over $\opSymm_+(2,\Omega)$.
\par
We now introduce the $(1,0)$-energy density and the $(1,0)$-energy as functionals defined over $\opGL_+(2,\Omega)$ and $\opSymm_+(2,\Omega)$. First, given an element $A\in\opGL_+(2,\Omega)$, its $(1,0)$-{\sl energy density} is defined by
$$
\hat{\Cal{E}}_g[A] = \sigma_g[A]\opdArea_g,\eqnum{\nexteqnno[DefnOneZeroEnergyDensity]}
$$
where $\opdArea_g$ here denotes the area form of $g$ and $\sigma_g(A)$ is the $(1,0)$-seminorm of $A$ with respect to $g$ as defined in Section \subheadref{LinearAlgebraicPreliminaries}. The $(1,0)$-{\sl energy} of any such $A$ is now defined by
$$
\hat{\opE}_g[A] := \int_\Omega\hat{\Cal{E}}_g[A].\eqnum{\nexteqnno[DefnTotalOneZeroEnergy]}
$$
By \eqnref{InvarianceOfSigma}, for all $M\in\opSO(2,\Omega)$ and for all $A\in\opGL_+(2,\Omega)$,
$$
\hat{\Cal{E}}_g[MA] = \hat{\Cal{E}}_g[A],
$$
so that both $\hat{\Cal{E}}_g$ and $\hat{\opE}_g$ descend to functionals defined over $\opSymm_+(2,\Omega)$. These functionals will also be referred to respectively as the $(1,0)$-{\sl energy density} and the $(1,0)$-{\sl energy} and will be denoted respectively by $\Cal{E}_g$ and $\opE_g$. In terms of the notation of \cite{TrapaniValli}, for all $h$,
$$
\opE_g[h] = T_\partial(g,h,\opId).\eqnum{\nexteqnno[ComparisonWithTV]}
$$
\newsubhead{Orbits of the diffeomorphism group}[OrbitsOfTheDiffeomorphismGroup]
Let $\chi_0(\Omega)$ denote the Frechet space of smooth vector fields $X:\Omega\rightarrow\Bbb{R}^2$ with compact support in $\Omega$. Let $\opDiff_0(\Omega)$ denote the group of smooth diffeomorphisms $\Phi:\Omega\rightarrow\Omega$ which coincide with the identity outside of some compact set. $\opDiff_0(\Omega)$ is a smooth Frechet Lie group locally modelled on $\chi_0(\Omega)$ whose Lie algebra also identifies with $\chi_0(\Omega)$. The pull-back operation defines a natural action of this group over $\opSymm_+(2,\Omega)$ and the corresponding infinitesimal action of $\chi_0(\Omega)$ on $\opSymm_+(2,\Omega)$ is given by the Lie derivative.
\par
We will henceforth be concerned with elements of $\opSymm_+(2,\Omega)$ which are critical for $\opE_g$ along the $\opDiff_0(\Omega)$ orbit in which they lie. Indeed, these are precisely the metrics $h$ for which the identity defines a {\sl one-harmonic map} from $(\Omega,g)$ to $(\Omega,h)$ in the sense of \cite{TrapaniValli}. A useful equivalent formulation is given as follows. Given a metric $h\in\opSymm_+(2,\Omega)$, the functional $\opE_{g,h}:\opDiff_0(\Omega)\rightarrow\Bbb{R}$ is defined by
$$
\opE_{g,h}[\Phi] := \opE_g[\Phi^* h].\eqnum{\nexteqnno[DefnEnergyOverOrbit]}
$$
The $L^2$-{\sl gradient} of $\opE_{g,h}$ at the identity, when it exists, is then defined to be the unique vector field $\nabla\opE_g[h]$ such that, for every other vector field $X\in\chi_0(\Omega)$,
$$
\opD\opE_{g,h}[\opId]\cdot X = \int_\Omega\langle\nabla\opE_g[h],X\rangle\opdArea_g.\eqnum{\nexteqnno[DefnLTwoGradient]}
$$
We will see presently that $\nabla\opE_g[h]$ exists for all $h$. In particular, $h$ is critical for $\opE_g$ along its $\opDiff_0(\Omega)$-orbit if and only if $\nabla\opE_g[h]$ vanishes.
\par
As mentioned above, it is easier to study functions over $\opSymm_+(2,\Omega)$ using its symmetric space structure. First, we readily show that, for all $A\in\opGL_+(2,\Omega)$ and for all $X\in\chi_0(\Omega)$, the Lie derivative of the metric $\oph[A]$ along the vector field $X$ is related to the derivative of the functional $\oph$ at $A$ by
$$
\Cal{L}_X\oph[A] = \opD\oph[A]\cdot\delta_X A,\eqnum{\nexteqnno[FormulaForLieDerivative]}
$$
where
$$
\delta_XA := \nabla_X A + A\nabla X.\eqnum{\nexteqnno[DefnDeltaXA]}
$$
\proclaim{Lemma \nextprocno}
\noindent For all $X$ and for all $A$, the derivatives of $\opE_{g,\oph[A]}$ and $\hat{\opE}_g$ are related by
$$
\opD\opE_{g,\oph[A]}[\opId]\cdot X = \opD\hat{\opE}_g[A]\cdot\delta_X A.\eqnum{\nexteqnno[VariationOfEInTermsOfDeltaA]}
$$
\endproclaim
\proof Indeed, by \eqnref{FormulaForLieDerivative},
$$\eqalign{
\opD\opE_{g,\oph[A]}[\opId]\cdot X &= \opD\opE_g[\oph[A]]\cdot\Cal{L}_X\oph[A]\cr
&= \opD\opE_g[\oph[A]]\cdot\opD\oph[A]\cdot \delta_X A\cr
&= \opD\hat{\opE}_g[A]\cdot\delta_X A,}
$$
as desired.\qed
\medskip
\noindent Using Stokes' theorem and \eqnref{FirstDerivativeA}, we obtain
\proclaim{Lemma \nextprocno}
\noindent For all $h\in\opSymm_+(2,\Omega)$, $\nabla\opE_g[h]$ exists. Furthermore, if $h=\oph[A]$, where $A$ has non-vanishing trace, then
$$
\nabla E_g[h] = - J\nabla\cdot(A J) + A J\nabla\frac{\opTr(AJ)}{\opTr(A)} + F(g,A)(\opTr(AJ),\nabla\opTr(AJ),\nabla\cdot(AJ)),\eqnum{\nexteqnno[DerivativeOfEnergy]}
$$
where $F(g,A)$ is a homogeneous quadratic polynomial which depends smoothly on $g$ and $A$. In particular, if $A$ is symmetric, then
$$
\nabla E_g[h] = - J\nabla\cdot(A J).\eqnum{\nexteqnno[DerivativeOfEnergyAtSymm]}
$$
\endproclaim
\proclabel{DerivativeOfEnergyAtSym}
\noindent This yields the following key fact.
\proclaim{Corollary \nextprocno}
\noindent An element $h\in\opSymm_+(2,\Omega)$ is a critical point of $\opE_g$ along its $\opDiff_0(\Omega)$-orbit if and only if
$$
\nabla\cdot(AJ)=0,\eqnum{\nexteqnno[CriticalPointCondition]}
$$
where here $A:=\opA[h]$. In particular, $h$ is a critical point of $\opE_g$ along its $\opDiff_0(\Omega)$-orbit if and only if $\opA[h]$ is a Codazzi field.
\endproclaim
\proclabel{CorDerivativeOfEnergyAtSym}
\noindent A lengthy, but straightforward, calculation likewise yields the first variation of $\nabla E_g[h]$ over $\opDiff_0(\Omega)$ orbits in $\opSymm_+(2,\Omega)$.
\proclaim{Lemma \nextprocno}
\noindent If $h\in\opSymm_+(2,\Omega)$ is a critical point of $E_g$ along its $\opDiff_0(\Omega)$ orbit then, for every vector field $X$ of compact support,
$$
\int_\Omega\langle\opD\nabla\opE_g[h]\cdot\Cal{L}_X h,X\rangle\opdArea_g
=\int_\Omega\frac{1}{\opTr(A)}\opTr(A(\nabla X) J)^2 - \kappa_g\langle X,AX\rangle\opdArea_g,\eqnum{\nexteqnno[SecondVariationFormula]}
$$
where here $A:=\opA[h]$ and $\kappa_g$ denotes the scalar curvature of $g$. In particular, when $g$ is negatively curved, this bilinear form is positive-definite.
\endproclaim
\proclabel{SecondVariationFormula}
\newsubhead{Recovering ellipticity}[RecoveringEllipticity]
By polarising \eqnref{SecondVariationFormula} and applying Stokes' Theorem to its right hand side, we see that, for all $g$ and $h$, the operator
$$
X\mapsto\opD\nabla\opE_g[h]\cdot\Cal{L}_Xh\eqnum{\nexteqnno[FirstComponentOfOperatorWRTX]}
$$
is a second order, linear partial differential operator. If this operator were elliptic, then \eqnref{SecondVariationFormula} would allow us to develop a local deformation theory for metrics which are critical points of $\opE_g$ along their $\opDiff_0(\Omega)$-orbits. The reality is, however, less simple. We now review the subtle idea used by Trapani \& Valli to recover ellipticity. First, for all $h\in\opSymm_+(2,\Omega)$, let $\kappa[h]$ denote its scalar curvature. By classical surface geometry, using, for example, the technique of moving frames, we readily obtain
\proclaim{Lemma \nextprocno}
\noindent If $h=\oph[A]$, then
$$
\opDet(A)\kappa[h] = \kappa[g] + \nabla\cdot JA^t\opDet(A)^{-1}\nabla\cdot(AJ).\eqnum{\nexteqnno[DifferentialCurvatureFormulaA]}
$$
In particular, when $A$ is symmetric
$$
\opDet(A)\kappa[h] = \kappa[g] + \nabla\cdot A^{-1}J\nabla\cdot(AJ).\eqnum{\nexteqnno[DifferentialCurvatureFormulaB]}
$$
\endproclaim
\proclabel{DifferentialCurvatureFormula}
\noindent For all $h\in\opSymm_+(\Omega)$, let $\opF_g[h]$ be the smooth function over $\Omega$ given by
$$
\opF_g[h] := -\frac{\kappa[h]\opdArea[h]}{\opdArea_g} + \kappa_g,\eqnum{\nexteqnno[DefinitionOfF]}
$$
and let $\opG_g[h]$ be the smooth vector field given by
$$
\opG_g[h] := \nabla\opF_g[h],\eqnum{\nexteqnno[DefinitionOfG]}
$$
where here the symbol $\nabla$ denotes the gradient with respect to the metric $g$, and is not to be confused with the $L^2$ gradient of the functional $\opF_g$. By Lemma \procref{DifferentialCurvatureFormula}, for all $h$, the vector field $\opG_g[h]$ is given explicitely by
$$
\opG_g[h] = -\nabla\opDet(A)^{-1}\nabla\cdot A^{-1} J\nabla\cdot(A J),\eqnum{\nexteqnno[ExplicitFormulaForF]}
$$
where here $A:=\opA[h]$. By diffeomorphism invariance of the curvature, for every diffeomorphism $\Phi$, we also have
$$
\opF_g[\Phi^*h]=-(\kappa[h]\opDet(A))\circ\Phi\frac{\Phi^*\opdArea_g}{\opdArea_g} + \kappa_g,\eqnum{\nexteqnno[DiffInvarianceOfF]}
$$
from which we readily obtain
\proclaim{Lemma \nextprocno}
\noindent If $h\in\opSymm_+(2,\Omega)$ is a critical point of $\opE_g$ along its $\opDiff_0(\Omega)$-orbit then, for all vector fields $X$ and $Y$ of compact support,
$$
\int_\Omega\langle\opD\opG_g[h]\cdot\Cal{L}_Xh,Y\rangle\opdArea_g
=\int_\Omega(\nabla\cdot(\kappa_gX))(\nabla\cdot Y)\opdArea_g.\eqnum{\nexteqnno[DerivativeOfG]}
$$
\endproclaim
\proclabel{DerivativeOfG}
\noindent It is now straightforward to show
\proclaim{Lemma \nextprocno}
\noindent If $\kappa_g<0$ and if $h\in\opSymm_+(2,\Omega)$ is a critical point of $\opE_g$ along its $\opDiff_0(\Omega)$ orbit, then the operator
$$
X\mapsto\opD(\nabla\opE_g - \opG_g)[h]\cdot\Cal{L}_Xh\eqnum{\nexteqnno[EllipticOperatorWRTX]}
$$
is a second-order, elliptic partial differential operator with self-adjoint principal symbol.
\endproclaim
\proclabel{EllipticityOfOperatorWRTX}
\remark[EllipticityOfOperatorWRTX] The alert reader might be suprised to find that \eqnref{EllipticOperatorWRTX} is of second order, since the formula for $\opG_g$ already involves three orders of differentiation so that \eqnref{EllipticOperatorWRTX} ought to be of fourth order. It is the diffeomorphism invariance of $\opF_g$ which, through \eqnref{DiffInvarianceOfF}, ensures that this operator is of second order.
\newsubhead{Perturbing one-harmonic diffeomorphisms}[PerturbingOneHarmonicDiffeomorphisms]
Let $\Surface$ be a compact Riemann surface of hyperbolic type with negatively-curved metric $g:=\langle\cdot,\cdot\rangle$. Let $\opSymm_+(2,\Surface)$ denote the space of smooth metrics over $\Surface$, let $\opGL_+(2,\Surface)$ denote the group of smooth, orientation-preserving endomorphism fields over $\Surface$, and let $\opSO(2,\Surface)$ denote its subgroup of endomorphism fields which are special orthogonal with respect to $g$. As before $\opSymm_+(2,\Surface)$ identifies with the symmetric space $\opSO(2,\Surface)\backslash\opGL_+(2,\Surface)$. Define $\oph:\opGL_+(2,\Surface)\rightarrow\opSymm_+(2,\Surface)$ by
$$
\oph[A] := g(A\cdot,A\cdot),
$$
and define the section $\opA:\opSymm_+(2,\Surface)\rightarrow\opGL_+(2,\Surface)$ as before. The functionals $\opE_g$, $\nabla\opE_g$, $\opF_g$ and $\opG_g$ are extended to functionals defined over $\opSymm_+(\Surface)$ in the usual manner using coordinate charts and partitions of unity. Finally, let $\chi(\Surface)$ denote the space of smooth vector fields over $\Surface$ and let $\opDiff(\Surface)$ denote the group of smooth diffeomorphisms of $\Surface$. As before, $\opDiff(\Surface)$ is a smooth Frechet Lie group locally modelled on $\chi(\Surface)$ with Lie algebra $\chi(\Surface)$.
\par
By \eqnref{DerivativeOfEnergyAtSymm}, \eqnref{ExplicitFormulaForF} and Stokes' theorem, we obtain
\proclaim{Lemma \nextprocno}
\noindent For all $h$, with $A:=\opA[h]$,
$$
\int_\Surface\langle\nabla\opE_g[h]-\opG_g[h],A^{-1}\nabla\opE_g[h]\rangle\opdArea_g
\geq \int_\Surface\langle\nabla\opE_g[h],A^{-1}\nabla\opE_g[h]\rangle\opdArea_g.\eqnum{\nexteqnno[ControlOfModifiedFunctional]}
$$
\endproclaim
\proclabel{ControlOfModifiedFunctional}
\noindent This yields
\proclaim{Lemma \nextprocno}
\noindent For all $h$, $\nabla\opE_g[h]$ vanishes if and only if $\nabla\opE_g[h]-\opG_g[h]$ vanishes. In particular, $h$ is a critical point of $\opE_g$ along its $\opDiff(\Surface)$-orbit if and only if $\nabla\opE_g[h]-\opG_g[h]$ vanishes. Furthermore, if $g$ is negatively curved then, for any such $h$, the operator
$$
X\mapsto D(\nabla\opE_g[h]-\opG_g[h])\cdot\Cal{L}_Xh\eqnum{\nexteqnno[PerturbationOperator]}
$$
has trivial kernel.
\endproclaim
\proclabel{EquivalenceOfModifiedFunctional}
\proof The only subtlety lies the final assertion. To prove this, we consider first the more general case of three vector-valued functions $u$, $v$ and $w$ which vanish at some point $x_0$ and which satisfy
$$
\langle u,v\rangle \geq \|w\|^2.
$$
Denoting their respective derivatives by $u'$, $v'$ and $w'$ and twice differentiating at $x_0$, we obtain
$$
\langle u',v'\rangle \geq \|w'\|^2.
$$
It follows that if $u'$ vanishes, then so too does $w'$. In particular, if $X$ is an element of the kernel of \eqnref{PerturbationOperator}, then
$$
D\nabla\opE_g[h]\cdot\Cal{L}_Xh = 0,
$$
so that, by \eqnref{SecondVariationFormula}, $X$ vanishes, as desired.\qed
\medskip
\noindent It follows that perturbations of zeros of $\nabla\opE_g$ are obtained by perturbing zeros $\nabla\opE_g-\opG_g$. This yields
\proclaim{Theorem \nextprocno}
\noindent Let $(g_t)_{t\in]-\epsilon,\epsilon[}$ and $(h_t)_{t\in]-\epsilon,\epsilon[}$ be smooth families of negatively curved metrics over $\Surface$. Let $\Phi:\Surface\rightarrow\Surface$ be a one-harmonic diffeomorphism from $(\Surface,g_0)$ to $(\Surface,h_0)$. Upon reducing $\epsilon$ if necessary, there exists a unique smooth family $(\Phi_t)_{t\in]-\epsilon,\epsilon[}$ such that $\Phi_0=\Phi$ and, for all $t$, $\Phi_t$ is a one-harmonic diffeomorphism from $(\Surface,g_t)$ to $(\Surface,h_t)$.
\endproclaim
\proclabel{PerturbationOfOneHarmonicDiffeomorphisms}
\proof Upon replacing $h_t$ with $\Phi^*h_t$ for all $t$, we may suppose that $\Phi=\opId$. For all $(k,\alpha)$, let $\chi^{k,\alpha}(\Surface)$ and $\opDiff^{k,\alpha}(\Surface)$ denote respectively the space of $C^{k,\alpha}$ vector fields over $\Surface$ and the space of $C^{k,\alpha}$ diffeomorphisms of $\Surface$. For a suitable neighbourhood $\Cal{U}^{k+2}$ of $0$ in $\chi^{k+2,\alpha}(\Surface)$, define $\Cal{E}:\Cal{U}^{k+2}\rightarrow\opDiff^{k+2,\alpha}(\Surface)$ by
$$
\Cal{E}[X](x) := \opExp(X(x)),
$$
where $\opExp$ here denotes the exponential map of some fixed metric over $\Surface$. Define the smooth functional $\Cal{F}:]-\epsilon,\epsilon[\times\Cal{U}^{k+2,\alpha}\rightarrow\chi^{k,\alpha}(\Surface)$ by
$$
\Cal{F}[t,X] := \nabla\opE_{g_t}[\Cal{E}(X)^*h_t] - \opG_{g_t}[\Cal{E}(X)^*h_t].
$$
The partial derivative at $(0,0)$ of this functional with respect to the second component is
$$
\opD_2\Cal{F}[0,0]\cdot X = \opD(\nabla\opE_{g_0} - \opG_{g_0})[h_0]\cdot\Cal{L}_Xh_0.
$$
By Lemma \procref{EllipticityOfOperatorWRTX}, this operator is elliptic with self-adjoint principal symbol and is therefore Fredholm of index zero. By Lemma \procref{EquivalenceOfModifiedFunctional}, it has trivial kernel and is thus invertible. It follows by the implicit function theorem that, upon reducing $\epsilon$ if necessary, there exists a unique smooth family $(X_t)_{t\in]-\epsilon,\epsilon[}\in\Cal{U}$ such that, for all $t$,
$$
\Cal{F}[t,X_t] = 0.
$$
For all $t$, denote $\Phi_t:=\Cal{E}(X_t)$. Smoothness of $\Phi_t$ follows by elliptic regularity, and this completes the proof.\qed
\newhead{Compactness}[Compactness]
\newsubhead{Controlling the conformal class}[ControllingTheConformalClass]
Let $\Surface$ be a compact surface of hyperbolic type furnished with an orientation form $\omega$. For the purposes of this chapter, a {\sl complex structure} over $\Surface$ is defined to be a smooth endomorphism field $J$ such that $J^2=-\opId$ and such that $\omega(\cdot,J\cdot)$ is positive definite. Let $\optildeConf(\Surface)$ denote the space of complex structures over $\Surface$ furnished with the $C^\infty$ topology. Let $\opConf(\Surface)$ denote its quotient space under the pull-back action of the diffeomorphism group $\opDiff(\Surface)$. In what follows, every complex structure in $\optildeConf(\Surface)$ will be identified with its equivalence class in $\opConf(\Surface)$.
\par
Recall by the Newlander-Nirenberg Theorem \cite{NewNir} that every two-dimensional complex structure is integrable so that $\opConf(\Surface)$ identifies with the moduli space of Riemann surfaces homeomorphic to $\Surface$. We emphasize that we are concerned in this chapter with {\sl unmarked} complex structures and that $\opConf(\Surface)$ is not to be confused with the Teichm\"uller space of {\sl marked} complex structures over $\Surface$.
\par
Let $\optildeHyp(\Surface)$ denote the space of hyperbolic metrics over $\Surface$ furnished with the $C^\infty$ topology. Let $\opHyp(\Surface)$ denote its quotient space under the pull-back action of $\opDiff(\Surface)$. Observe that the topology of $\opHyp(\Surface)$ is precisely the $C^\infty$ Cheeger-Gromov topology. In what follows, every metric in $\optildeHyp(\Surface)$ will also be identified with its equivalence class in $\opHyp(\Surface)$.
\par
Recall that $\opConf(\Surface)$ and $\opHyp(\Surface)$ are canonically homeomorphic. Indeed, a complex structure $J$ and a metric $h$ are said to be {\sl compatible} with one another whenever
$$
h(J\cdot,J\cdot) = h(\cdot,\cdot).
$$
There is a unique complex structure $J$ compatible with every metric $h$. Conversely, by Riemann's uniformisation theorem, there is a unique hyperbolic metric $h$ compatible with every complex structure $J$. This yields an identification of $\optildeHyp(\Surface)$ with $\optildeConf(\Surface)$ which descends to the desired homeomorphism.
\par
Consider now a sequence $(g_m,A_m)_{m\in\Bbb{N}}$ where, for all $m$, $g_m$ is a riemannian metric over $\Surface$ and $A_m$ is a positive-definite endomorphism field which is symmetric with respect to $g_m$. For all $m$, let $J_m$ denote the complex structure of $g_m$ and define
$$
\hat{J}_m := \frac{1}{\sqrt{\opDet(A_m)}}J_mA_m.\eqnum{\nexteqnno[DefinitionOfIntermediateComplexStructure]}
$$
It follows from \eqnref{TwoDimMatrixRelations} that, for all $m$, $\hat{J}_m$ is a complex structure which is compatible with the metric
$$
h_m := g_m(A_m\cdot,\cdot).\eqnum{\nexteqnno[DefinitionOfIntermediateMetric]}
$$
In this section we establish conditions under which the sequence $(\hat{J}_m)_{m\in\Bbb{N}}$ subconverges in $\opConf(\Surface)$. First, for any metric $h$ over $\Surface$, define its {\sl systole} by
$$
\opSys(h) := \minf_\gamma\opLength(\gamma,h),\eqnum{\nexteqnno[DefinitionOfSystole]}
$$
where the infimum is taken over all homotopically non-trivial closed curves $\gamma$ in $\Surface$ and, for all such $\gamma$, $\opLength(\gamma,h)$ denotes its length with respect to $h$. Observe that, when $h$ is non-positively curved, the systole is equal to twice the injectivity radius. In this section, we prove
\proclaim{Lemma \nextprocno}
\noindent If there exists $B>0$ such that, for all $m$,
$$\eqalign{
\opSys(g_m) &\geq \frac{1}{B},\cr
\opDet(A_m) &\geq \frac{1}{B}\ \text{and}\cr
\int_\Surface\opTr(A_m)\opdArea_m &\leq B\vphantom{\frac{1}{B}},\cr}
$$
then the sequence $(\hat{J}_m)$ is precompact in  $\opConf(\Surface)$.
\endproclaim
\proclabel{CompactnessOfIntermediateConformalClasses}
We first recall the following precompactness criterion for hyberbolic metrics (see \cite{Mumford}).
\proclaim{Theorem \nextprocno, {\bf Mumford (1971)}}
\noindent A subset $X$ of $\opHyp(\Surface)$ is precompact if and only if
$$
\minf_{h\in X}\opSys(h) > 0.
$$
\endproclaim
\proclabel{PrecompactnessForHyperbolicMetrics}
\noindent In order to prove Lemma \procref{CompactnessOfIntermediateConformalClasses}, it suffices to express Theorem \procref{PrecompactnessForHyperbolicMetrics} in conformal terms. To this end, we recall the concept of conformal modulus. First, a {\sl conformal annulus} is defined to be a Riemann surface $S$ with fundamental group $\Bbb{Z}$. Up to conformal equivalence, the only conformal annulus of parabolic type is $\Bbb{C}^*$. Every other conformal annulus has hyperbolic type and is conformally equivalent to $S^1\times]0,M[$ for some unique $M\in]0,\infty]$. The quantity $M$, which may be considered as a conformal ``length'' of the annulus $S$, is called the {\sl conformal modulus} of $S$ and will be denoted by $\opMod(S)$. It is also determined by the following formula (c.f. \cite{Ahlfors}),
$$
\frac{1}{\opMod(S)} = \msup_g\frac{\opSys(g)^2}{2\pi\opArea(g)},\eqnum{\nexteqnno[ExtremalLengthFormula]}
$$
where the supremum is taken over all metrics $g$ over $S$ which are compatible with the complex structure.
\par
An open subset $\Omega$ of $\Surface$ will be said to be an {\sl essential annulus} whenever it has fundamental group $\Bbb{Z}$ and every non-contractible curve in $\Omega$ is also non-contractible in $\Surface$. Given a complex structure $J$ over $\Surface$, we define
$$
\opMod(J) := \msup_{\Omega\subseteq\Surface}\opMod(\Omega,J),\eqnum{\nexteqnno[DefinitionOfModuleOfConformalStructure]}
$$
where the supremum is taken over all essential annuli $\Omega$ in $\Surface$ and, for all such $\Omega$, $(\Omega,J)$ denotes the Riemann surface defined by restricting $J$ to $\Omega$. Heuristically, $\opMod(J)$ is the supremal ``length'' with respect to $J$ of essential annuli in $\Surface$.
\par
From \eqnref{ExtremalLengthFormula} and the Gauss-Bonnet theorem, we obtain
\proclaim{Lemma \nextprocno}
\noindent If $J$ is a complex structure over $\Surface$ with hyperbolic metric $h$, then
$$
\opMod(J) \leq \frac{4\pi^2\left|\chi(\Surface)\right|}{\opSys(h)^2},\eqnum{\nexteqnno[UpperBoundOfConformalModule]}
$$
where $\chi(\Surface)$ here denotes the Euler characteristic of $\Surface$.
\endproclaim
\proclabel{UpperBoundOfConformalModule}
\remark[UpperBoundOfConformalModule] We are not aware of whether this bound is sharp.
\medskip
\noindent Conversely,
\proclaim{Lemma \nextprocno}
\noindent There exists a function $F:]0,\infty[\rightarrow]0,\infty[$ such that, if $J$ is a complex structure over $\Surface$ with hyperbolic metric $h$, then
$$
\opMod(J) \geq F(\opSys(h)).\eqnum{\nexteqnno[LowerBoundOfConformalModule]}
$$
Furthermore, $F(x)$ tends to infinity as $x$ tends to zero.
\endproclaim
\proclabel{LowerBoundOfConformalModule}
\proof Let $J$ be a complex structure over $\Surface$ with hyperbolic metric $h$. By Theorem $4.1.1$ of \cite{Buser}, if $\gamma$ is a simple, closed geodesic of $h$ of length $2l_1$, then distinct lifts of $\gamma$ in $\Bbb{H}^2$ are separated by at least $l_2$, where
$$
\opSinh(l_1)\opSinh(l_2) \geq 1.\eqnum{\nexteqnno[UncertaintyPrinciple]}
$$
Denote $\epsilon:=\opSys(h)$. Let $\gamma$ be a simple, closed geodesic of $h$ of length $\epsilon$. Choose $R>0$ such that
$$
\opSinh(\epsilon/2)\opSinh(2R) < 1.
$$
By \eqnref{UncertaintyPrinciple}, the geodesic segments of length $2R$ which meet $\gamma$ orthogonally at their midpoints foliate an open subset $\Omega$ of $\Surface$. This subset is an essential annulus which, by classical hyperbolic geometry, is isometric to the cylinder $S^1\times ]-R,R[$ furnished with the metric
$$
\frac{\epsilon^2\opCosh^2(r)}{4\pi^2}d\theta^2 + dr^2.
$$
This is in turn conformally equivalent to the cylinder $S^1\times]0,2L[$ furnished with the product metric, where
$$
L := \frac{2\pi}{\epsilon}\opArcTan(\opTanh(R)).
$$
The result now follows.\qed
\medskip
\noindent Theorem \procref{PrecompactnessForHyperbolicMetrics} and Lemmas \procref{UpperBoundOfConformalModule} and \procref{LowerBoundOfConformalModule} together yield
\proclaim{Theorem \nextprocno}
\noindent A subset $X$ of $\opConf(\Surface)$ is precompact if and only if
$$
\msup_{J\in X}\opMod(J) < \infty.
$$
\endproclaim
\proclabel{PrecompactnessForComplexStructures}
\noindent We now prove Lemma \procref{CompactnessOfIntermediateConformalClasses}.
\medskip
{\bf\noindent Proof of Lemma \procref{CompactnessOfIntermediateConformalClasses}:\ }For all $m$, define
$$
\hat{h}_m := \frac{\opTr(A_m)}{\sqrt{\opDet(A_m)}}h_m.
$$
For all $m$, the hypotheses on $(g_m)$ and $(A_m)$ ensure that,
$$\eqalign{
\opSys(\hat{h}_m) &\geq B^{-5/4}\ \text{and}\cr
\opArea(\hat{h}_m) &\leq B,\cr}
$$
so that, by \eqnref{ExtremalLengthFormula},
$$
\opMod(\hat{J}_m) \leq 2\pi B^4.
$$
The result now follows by Theorem \procref{PrecompactnessForComplexStructures}.\qed
\newsubhead{Harmonic maps}[HarmonicMaps]
We now recall the definition of harmonicity for smooth maps between riemannian manifolds. We refer the reader to \cite{Jost} for a complete treatment of this theory. Let $M$ and $N$ be riemannian manifolds with respective metrics $g$ and $h$. Let $\Phi:M\rightarrow N$ be a smooth function. The derivative $D\Phi$ of $\Phi$ defines a smooth section of $\opLin(TM,\Phi^*TN)$ over $M$. The Levi-Civita covariant derivatives of $g$ and $h$ define a covariant derivative over this bundle. The {\sl laplacian} of $\Phi$ is then defined to be the section of $\Phi^*TN$ given by
$$
\Delta\Phi := \sum_{i=1}^m\big(\nabla_{e_i}D\Phi)e_i,\eqnum{\nexteqnno[DefinitionOfLaplacian]}
$$
where $e_1,...,e_m$ is a local orthonormal frame of $M$. We readily verify that \eqnref{DefinitionOfLaplacian} is independent of the orthonormal frame chosen and is thus well-defined. In coordinate charts, $\Delta\Phi$ is given by (see $(8.1.16)$ of \cite{Jost})
$$
(\Delta\Phi)^k = g^{mn}\Phi^k_{;mn} - g^{mn}\gamma^p_{mn}\Phi^k_{;p} + g^{mn}(\Gamma^k_{pq}\circ\Phi)\Phi^p_{;m}\Phi^q_{;n},\eqnum{\nexteqnno[LaplacianInCoordinateCharts]}
$$
where $\gamma$ and $\Gamma$ here denote respectively the Christoffel symbols of $g$ and $h$ and the summation convention is assumed. The function $\Phi$ is said to be {\sl harmonic} whenever
$$
\Delta\Phi = 0.\eqnum{\nexteqnno[DefinitionOfHarmonicity]}
$$
In what follows, we will be concerned with the case where $M$ is a surface. Here, it is known (see \cite{Jost}) that if two metrics $g$ and $g'$ are conformally equivalent to one another, then $\Phi$ is harmonic with respect to the one if and only if it is harmonic with respect to the other. In other words, the property of harmonicity only depends on the complex structure of $g$.
\par
The compactness results we aim to prove will require a generalisation of this concept involving a new term on the right-hand side. Although this generalisation no longer arises as the Euler-Lagrange equations of some functional, it is still conformally invariant in the case where $M$ is $2$-dimensional. It is this latter property which ensures that compactness results similar to those known for harmonic maps continue to hold.
\par
Suppose that $N$ is the cartesian product $N_1\times...\times N_k$ of finitely many riemannian manifolds. Let $\alpha:=(\alpha_1,...,\alpha_k)$ be a vector of $1$-forms over $N$. A smooth function $\Phi:M\rightarrow N$ will be said to be $\alpha$-harmonic whenever, for each $i$,
$$
\Delta\Phi_i = D\Phi_i\cdot(\Phi^*\alpha_i)^\#,\eqnum{\nexteqnno[FHarmonicity]}
$$
where $\Phi_1,...,\Phi_k$ denote the components of $\Phi$ and $\#$ denotes Berger's musical isomorphism with respect to the metric $g$. In coordinate charts, the right-hand side of \eqnref{FHarmonicity} is given by
$$
(D\Phi_i\cdot(\Phi^*\alpha_i)^\#)^k = g^{mn}(\alpha_{i,p}\circ\Phi)\Phi^p_{;m}\Phi^k_{i;n}.\eqnum{\nexteqnno[FHarmonicityInCoordinates]}
$$
In particular, when $M$ is two-dimensional, $\alpha$-harmonicity also only depends on the complex structure of $g$.
\par
We now prove a compactness result for $\alpha$-harmonic functions. Recall first that the {\sl energy} of a smooth function $\Phi:M\rightarrow N$ is defined by
$$
\opE(\Phi) := \int_M\|D\Phi\|^2\opdVol_M,\eqnum{\nexteqnno[DefinitionOfEnergyForFunctions]}
$$
where $\opdVol_M$ denotes the volume form of $M$ and, in every coordinate chart,
$$
\|D\Phi\|^2 := g^{mn}h_{pq}\Phi^p_m\Phi^q_n.\eqnum{\nexteqnno[DefinitionOfEnergyForFunctionsII]}
$$
When $M$ is two-dimensional, $\opE(\Phi)$ also only depends on the complex structure of $g$. We show
\proclaim{Theorem \nextprocno}
\noindent Let $\Surface$ be a compact riemannian surface. Let $N:=N_1\times...\times N_k$ be a cartesian product of compact, non-positively curved manifolds. Let $\alpha:=(\alpha_1,...,\alpha_k)$ be a vector of smooth $1$-forms over $N$. Let $(\Phi_m)_{m\in\Bbb{N}}$ be a sequence of $\alpha$-harmonic functions from $\Surface$ into $N$. If
$$
\msup_m \opE(\Phi_m) <\infty,\eqnum{\nexteqnno[CompactnessForFHarmonicFunctions]}
$$
then there exists an $\alpha$-harmonic function $\Phi_\infty$ towards which $(\Phi_m)$ subconverges in the $C^\infty$ sense.
\endproclaim
\proclabel{CompactnessForFHarmonicFunctions}
\remark[CompactnessForFHarmonicFunctionsI] For ease of presentation, Theorem \procref{CompactnessForFHarmonicFunctions} is stated in its simplest form. In its most general form, the metric of $\Surface$, the metrics of $N_1,...,N_k$ and the $1$-forms $\alpha_1,...,\alpha_k$ are allowed to vary within suitable compact sets. In addition, $\Surface$ and $N_1,...,N_k$ can be taken to be non-compact provided suitable conditions hold at infinity.
\medskip
Theorem \procref{CompactnessForFHarmonicFunctions} follows from the observation that, for any $\alpha$-harmonic function $\Phi:\Surface\rightarrow N$, for each $i$, and for any smooth, convex function $F:N_i\rightarrow\Bbb{R}$,
$$
\Delta(F\circ\Phi_i) - d(F\circ\Phi_i)((\Phi^*\alpha_i)^\#) \geq 0,\eqnum{\nexteqnno[PreMaximumPrinciple]}
$$
so that $(F\circ\Phi_i)$ satisfies a strong maximum principle (c.f. \cite{GilbTrud}). The remainder of the proof follows as in the case of harmonic functions. First, we obtain the following Liouville-type result.
\proclaim{Lemma \nextprocno}
\noindent Let $N:=N_1\times...\times N_k$ be a cartesian product of complete manifolds of non-positive curvature. Let $\alpha:=(\alpha_1,...,\alpha_k)$ be vector of smooth $1$-forms over $N$. There exists no non-constant $\alpha$-harmonic function $\Phi:\Bbb{R}^2\rightarrow N$ of finite energy.
\endproclaim
\proclabel{LiouvilleTheorem}
\proof Let $\Phi:\Bbb{R}^2\rightarrow N$ be an $\alpha$-harmonic function of finite energy. Upon lifting $\Phi$ to the universal cover, we may suppose that $N_i$ is simply connected for all $i$. For all $r>0$, let $S_r$ and $D_r$ denote respectively the circle and disk of radius $r$ about $0$ in $\Bbb{R}^2$. The Cauchy-Schwartz inequality yields
$$\eqalign{
\int_1^\infty r\opLength(\Phi(S_r))^2dr &\leq \int_1^\infty\frac{1}{r}\bigg(\int_0^{2\pi}\|D\Phi\|r d\theta\bigg)^2dr\cr
&\leq \int_1^\infty\int_0^{2\pi}2\pi\|D\Phi\|^2 rdrd\theta\cr
&\leq \opE(\Phi),\vphantom{\frac{1}{2}}\cr}
$$
so that, for all $i$
$$
\mliminf_{r\rightarrow\infty}\opDiam(\Phi_i(S_r)) \leq \mliminf_{m\rightarrow\infty}\opLength(\Phi(S_r)) = 0.
$$
However, for all $i$, since $N_i$ is simply-connected and non-positively curved, the squared distance to any point in $N_i$ is a smooth, convex function so that, by the maximum principle,
$$
\mliminf_{m\rightarrow\infty}\opDiam(\Phi_i(D_r)) = 0.
$$
The result follows.\qed
\medskip
\noindent A standard elliptic bootstrapping argument yields (see \cite{GilbTrud})
\proclaim{Theorem \nextprocno}
\noindent Let $\Surface$ be a riemannian surface. Let $N:=N_1\times...\times N_k$ be a cartesian product of compact riemannian manifolds. Let $\alpha:=(\alpha_1,...,\alpha_k)$ be a vector of smooth $1$-forms over $N$. Let $(\Phi_m)_{m\in\Bbb{N}}$ be a sequence of $\alpha$-harmonic functions from $\Surface$ into $N$. If
$$
\msup_m \|D\Phi_m\|_{L^\infty} < \infty,\eqnum{\nexteqnno[COneCompactness]}
$$
then there exists an $\alpha$-harmonic function $\Phi_\infty$ towards which $(\Phi_m)$ subconverges in the $C^\infty_\oploc$ sense.
\endproclaim
\proclabel{COneCompactness}
\remark[COneCompactness] As before, for ease of presentation, Theorem \procref{COneCompactness} is stated in its simplest form. In its most general form, the metric of $\Surface$, the metrics of $N_1,...,N_k$ and the $1$-forms $\alpha_1,...,\alpha_k$ are allowed to vary within suitable compact sets. In addition, if there exists a point $x\in\Surface$ such that $(\Phi_m(x))_{m\in\Bbb{N}}$ converges in $N$, then $N_1,...,N_k$ may also be taken to be non-compact.
\medskip
\noindent A blow-up argument now allows us to conclude.
\medskip
{\bf\noindent Proof of Theorem \procref{CompactnessForFHarmonicFunctions}:\ }In light of Theorem \procref{COneCompactness}, it suffices to show that
$$
\msup_m \|D\Phi_m\|_{L^\infty}<\infty.
$$
Suppose the contrary. We may suppose that there exists a sequence $(x_m)$ of points in $\Surface$ and a sequence $(R_m)$ of positive real numbers converging to infinity such that, for all $m$,
$$
\|D\Phi_m\|_{L^\infty} = \|D\Phi_m(x_m)\| = R_m.
$$
For ease of presentation, we may suppose that $x_m=x_0$ for all $m$. Furthermore, choosing a suitable chart of $\Surface$ about $x_0$, we may suppose that $\Surface$ is the unit disk in $\Bbb{R}^2$, $x_0=0$ and
$$
g_{ij} = e^{2\phi}\delta_{ij},
$$
where $\phi(0)=0$. For all $m$, define $\tilde{\Phi}_m:B_{R_m}(0)\rightarrow N$ by
$$
\tilde{\Phi}_m(x) := \Phi_m(x/R_m),
$$
and define the metric $g_m$ over $B_{R_m}(0)$ by
$$
g_m(x)_{ij} := e^{2\phi(x/R_m)}\delta_{ij}.
$$
Observe that $(g_m)$ converges smoothly over every compact subset of $\Bbb{R}^2$ to the euclidean metric. By Theorem \procref{COneCompactness} and the subsequent remark, $(\tilde{\Phi}_m)_{m\in\Bbb{N}}$ subconverges in the $C^\infty_\oploc$ sense to an $\alpha$-harmonic function $\tilde{\Phi}_\infty:\Bbb{R}^2\rightarrow N$. Since $D\tilde{\Phi}_\infty(0)$ has unit norm, this function is non-constant. However, by conformal invariance of the energy,
$$
\opE(\tilde{\Phi}_\infty) \leq \mlimsup_{m\rightarrow\infty}\opE(\Phi_m) < \infty,
$$
which is absurd by Lemma \procref{LiouvilleTheorem}. This completes the proof.\qed
\newsubhead{Compactness}[Compactness]
\noindent Let $\Surface$ be a compact surface. Let $(g_{1,m})_{m\in\Bbb{N}}$ and $(g_{2,m})_{m\in\Bbb{N}}$ be sequences of negatively curved metrics over $\Surface$. For all $m$, let $\Phi_m:(\Surface,g_{1,m})\rightarrow(\Surface,g_{2,m})$ be a one-harmonic diffeomorphism and let $A_m$ be the unique positive-definite Codazzi field such that
$$
\Phi_m^*g_{2,m} = g_{1,m}(A_m\cdot,A_m\cdot).\eqnum{\nexteqnno[CompactnessDefinitionOfAm]}
$$
For all $m$, let $\kappa_{1,m}$ and $\kappa_{2,m}$ denote the curvatures of $g_{1,m}$ and $g_{2,m}$ respectively, so that, by \eqnref{DifferentialCurvatureFormulaB},
$$
\kappa_{2,m}\circ\Phi_m = \kappa_{1,m} e^{-\phi_m},\eqnum{\nexteqnno[CompactnessBigPhiAndLittlePhi]}
$$
where
$$
\phi_m := \opLog(\opDet(A_m)).\eqnum{\nexteqnno[CompactnessDefnOfLittlePhi]}
$$
In this section, we first recover the compactness result of Trapani \& Valli, namely
\proclaim{Theorem \nextprocno}
\noindent If $(g_{1,m})_{m\in\Bbb{N}}$ and $(g_{2,m})_{m\in\Bbb{N}}$ converge smoothly to negatively curved metrics $g_{1,\infty}$ and $g_{2,\infty}$ respectively, and if
$$
\msup_m\opE(\Phi_m) < \infty,
$$
then there exists a smooth one-harmonic diffeomorphism $\Phi_\infty:(\Surface,g_{1,\infty})\rightarrow(\Surface,g_{2,\infty})$ towards which $(\Phi_m)_{m\in\Bbb{N}}$ subconverges smoothly.
\endproclaim
\proclabel{CompactnessTrapaniValli}
\noindent Using almost identical techniques, we also show
\proclaim{Theorem \nextprocno}
\noindent If $(g_{1,m})_{m\in\Bbb{N}}$ converges smoothly to the negatively curved metric $g_{1,\infty}$, if $(\phi_m)_{m\in\Bbb{N}}$ converges smoothly to the function $\phi_\infty$, and if
$$
\msup_m\int_\Surface\opTr(A_m)\opdArea_m < \infty,
$$
then there exists a symmetric, positive-definite Codazzi field $A_\infty$ of $g_{1,\infty}$ towards which $(A_m)_{m\in\Bbb{N}}$ subconverges smoothly.
\endproclaim
\proclabel{CompactnessForCodazziFields}
As in Section \subheadref{ControllingTheConformalClass}, for all $m$, define
$$
\hat{J}_m := e^{-\phi_m/2}J_mA_m,\eqnum{\nexteqnno[CompactnessIntermediateConformalStructure]}
$$
where $J_m$ denotes the complex structure of $g_m$. Observe that, by Lemma \procref{CompactnessOfIntermediateConformalClasses}, the hypotheses of both Theorems \procref{CompactnessTrapaniValli} and \procref{CompactnessForCodazziFields} imply that there exists a complex structure $\tilde{J}_\infty$ over $\Surface$ and a sequence $(\tilde{\Phi}_m)_{m\in\Bbb{N}}$ of smooth diffeomorphisms of $\Surface$ such that the sequence $(\tilde{\Phi}_m^*\hat{J}_m)_{m\in\Bbb{N}}$ converges smoothly to $\tilde{J}_\infty$. For all $m$, denote
$$
\tilde{J}_m := \tilde{\Phi}_m^*\hat{J}_m,\eqnum{\nexteqnno[CompactnessDefinitionOfTildeJ]}
$$
and define the $1$-form $\alpha_m$ over $\Surface$ by
$$
\alpha := -\frac{1}{2}d\phi_m.\eqnum{\nexteqnno[CompactnessDefinitionOfAlpha]}
$$
\proclaim{Lemma \nextprocno}
\noindent For all $m$, $\tilde{\Phi}_m$ is an $\alpha_m$-harmonic function from $(\Surface,\tilde{J}_m)$ into $(\Surface,g_{1,m})$.
\endproclaim
\proclabel{CompactnessAlphaHarmonicity}
\proof Observe that this is equivalent to proving that the identity is an $\alpha_m$-harmonic function from $(\Surface,\hat{J}_m)$ into $(\Surface,g_{1,m})$. Now, choose $m$ and denote $g:=g_m$, $\hat{J}:=\hat{J}_m$, $A:=A_m$ and $\alpha:=\alpha_m$. Define the metric $h$ by
$$
h := g(A\cdot,\cdot)
$$
and observe that $\hat{J}$ is the complex structure of $h$. Now choose a chart of $\Surface$. Since $A$ is a symmetric Codazzi field, the tensor
$$
A_{ijk} := g_{ip}A^p_{j;k}
$$
is symmetric in its three components. It follows by the Koszul formula that the difference between the Christoffel symbols of $h$ and $g$ is
$$
\Gamma^k_{ij} = \frac{1}{2}h^{kp}A_{pij}.
$$
The laplacian of the identity map from $(\Surface,h)$ to $(\Surface,g)$ is thus
$$
(\Delta\opId)^k = -\frac{1}{2}h^{ij}h^{kp}A_{pij}
= -\frac{1}{2}h^{kp}\opLog(\opDet(A))_{;p}
= h^{pq}(\alpha_r\circ\opId)\delta^r_p\delta^k_q,
$$
as desired.\qed
\medskip
It is now helpful to view the problem more symmetrically. Thus, for all $m$, define $\Psi_m:\Surface\rightarrow\Surface\times\Surface$ by
$$
\Psi_m(x) := (\tilde{\Phi}_m(x),(\Phi_m\circ\tilde{\Phi}_m)(x)).\eqnum{\nexteqnno[CompactnessDefinitionOfPsi]}
$$
For all $m$, define
$$
\eqalign{
\beta_{1,m} &:= \frac{1}{2}d\opLog(\kappa_{2,m}(y)) - \frac{1}{2}d\opLog(\kappa_{1,m}(x))\ \text{and}\cr
\beta_{2,m} &:= \frac{1}{2}d\opLog(\kappa_{1,m}(x)) - \frac{1}{2}d\opLog(\kappa_{2,m}(y)),\cr}\eqnum{\nexteqnno[CompactnessDefinitionOfBetaI]}
$$
and define
$$
\beta_m := (\beta_{1,m},\beta_{2,m}).\eqnum{\nexteqnno[CompactnessDefinitionOfBetaII]}
$$
For all $m$,
$$
\Psi_m^*\beta_{1,m} = \tilde{\Phi}_m^*\alpha_m,
$$
so that, by Lemma \procref{CompactnessAlphaHarmonicity},
$$
\Delta\Psi_{1,m} = D\Psi_{1,m}\cdot(\Psi_m^*\beta_{1,m}).\eqnum{\nexteqnno[CompactnessBetaHarmonicityI]}
$$
Likewise, by symmetry,
$$
\Delta\Psi_{2,m} = D\Psi_{2,m}\cdot(\Psi_m^*\beta_{2,m}).\eqnum{\nexteqnno[CompactnessBetaHarmonicityII]}
$$
It follows that, for all $m$, $\Psi_m$ is a $\beta_m$-harmonic function from $(\Surface,\tilde{J}_m)$ into $(\Surface\times\Surface,g_{1,m}\oplus g_{2,m})$. Furthermore, a straightforward calculation yields, for all $m$,
$$\eqalign{
\opE(\Psi_{1,m}) &= \int_\Surface\opTr(A_m)e^{-\frac{1}{2}\phi_m}\opdVol_m\ \text{and}\cr
\opE(\Psi_{2,m}) &= \int_\Surface\opTr(A_m)e^{\frac{1}{2}\phi_m}\opdVol_m,\cr}\eqnum{\nexteqnno[CompactnessEnergyBoundsI]}
$$
so that
$$
\opE(\Psi_m) = \opE(\Psi_{1,m}) + \opE(\Psi_{2,m}) = \int_\Surface2\opTr(A_m)\opCosh(\phi_m/2)\opdVol_m.\eqnum{\nexteqnno[CompactnessEnergyBoundsII]}
$$
We now prove Theorem \procref{CompactnessTrapaniValli}.
\medskip
{\bf\noindent Proof of Theorem \procref{CompactnessTrapaniValli}:\ }By Theorem \procref{CompactnessForFHarmonicFunctions} and the subsequent remark, there exists a smooth function $\Psi_\infty:\Surface\rightarrow\Surface\times\Surface$ towards which $(\Psi_m)_{m\in\Bbb{N}}$ subconverges smoothly. Let $\Psi_{1,\infty}$ and $\Psi_{2,\infty}$ denote its two components. Choose $i\in\left\{1,2\right\}$ and consider the function
$$
f:=\opLog(\left|\partial\Psi_{i,\infty}\right|^2) - \opLog(\left|\overline{\partial}\Psi_{i,\infty}\right|^2).
$$
This function is non-negative and vanishes if and only if $D\Psi_{i,\infty}$ is degenerate. In \cite{TrapaniValli}, Trapani \& Valli use a maximum principle to show that if this function vanishes at a single point, then it vanishes identically. However, by compactness of $\Surface$, $\Psi_{i,\infty}$ is surjective so that $D\Psi_{i,\infty}$ is non-degenerate at at least one point. It follows that both $\Psi_{1,\infty}$ and $\Psi_{2,\infty}$ are smooth diffeomorphisms. In particular, $(\Phi_m)_{m\in\Bbb{N}}$ converges smoothly to $(\Psi_{2,\infty}\circ\Psi_{1,\infty}^{-1})$ and the result follows.\qed
\medskip
{\bf\noindent Proof of Theorem \procref{CompactnessForCodazziFields}:\ }By Theorem \procref{CompactnessForFHarmonicFunctions} and the subsequent remark, there exists a smooth function $\tilde{\Phi}_\infty$ towards which $(\tilde{\Phi}_m)_{m\in\Bbb{N}}$ subconverges smoothly. As in the proof of Theorem \procref{CompactnessTrapaniValli}, $\tilde{\Phi}_\infty$ is a smooth diffeomorphism. By \eqnref{CompactnessIntermediateConformalStructure}, $(A_m)_{m\in\Bbb{N}}$ subconverges smoothly to
$$
A_\infty := -e^{\phi_\infty/2}J_\infty((\tilde{\Phi}_\infty)_*\tilde{J}_\infty),
$$
where $J_\infty$ here denotes the complex structure of $g_\infty$. This completes the proof.\qed
\medskip
\noindent For completeness, we close this chapter by recalling the proof of Trapani \& Valli's result.
\proclaim{Theorem \nextprocno}
\noindent Let $g_1$ and $g_2$ be marked, negatively-curved metrics. There exists a unique one-harmonic diffeomorphism $\Phi:(S,g_1)\rightarrow (S,g_2)$ which preserves the marking.
\endproclaim
\proclabel{MainExistenceThmTP}
\proof We may suppose that the marking of $g_1$ coincides with that of $g_2$, so that we are interested in diffeomorphisms which are isotopic to the identity. For each $i$, let $u_i$ solve
$$
\Delta^i u_i = \kappa_i + e^{2u_i},
$$
where $\kappa_i$ and $\Delta^i$ denote respectively the curvature of $g_i$ and its Laplace-Beltrami operator. For all $t\in[0,1]$, define
$$
g_{i,t} := e^{2(1-t)u_i}g_i.
$$
Then, for each $i$, $g_{i,0}=g_i$, $g_{i,1}$ is hyperbolic and, for all $t$, $g_{i,t}$ is negatively curved. Let $I\subseteq[0,1]$ be the set of all $t$ for which there exists a unique one-harmonic diffeomorphism $\Phi_t:(\Surface,g_{1,t})\rightarrow(\Surface,g_{2,t})$ isotopic to the identity. Since $g_{1,1}$ and $g_{2,1}$ have constant curvature, by \cite{Schoen} and \cite{SchoenYau}, $1\in I$. By Theorems \procref{PerturbationOfOneHarmonicDiffeomorphisms} and \procref{CompactnessTrapaniValli}, $I$ is both open and closed. It follows by connectedness that $0\in I$, and existence and uniqueness follow.\qed
\newhead{The Weyl problem}[TheWeylProblem]
\newsubhead{The function $\hat{\opE}$ and its derivatives}[TheDerivativesOfE]
Let $\Surface$ be a compact surface of hyperbolic type. Given two marked, negatively-curved metrics $g$ and $h$ over $\Surface$, define
$$
\hat{\opE}[g,h] := \int_\Surface\opTr(A)\opdArea_g,\eqnum{\nexteqnno[TwoArgumentEnergyFunctional]}
$$
where $\opdArea_g$ here denotes the area form of $g$ and $A$ is the unique positive-definite Codazzi field given by Theorem \procref{TrapaniValliII}. Observe that, for all $g$ and $h$,
$$
\hat{\opE}[g,h] = \hat{\opE}[h,g].\eqnum{\nexteqnno[SymmetryOfTwoArgumentEnergyFunctional]}
$$
We now consider the case where $g$ is an arbitrary negatively-curved metric and $h$ is hyperbolic. We take $g$ to be fixed, and we define the function $\hat{\opE}_g$ over $\opTeich[\Surface]$ by
$$
\hat{\opE}_g[h] := \hat{\opE}[g,h].\eqnum{\nexteqnno[DefinitionOfFunctionOverTeichmuellerSpace]}
$$
In Sections $3$ and $5$ of \cite{BonMonSch}, Bonsante, Mondello \& Schlenker study the first and second derivatives of this function in the case where $g$ also has constant curvature. In particular, they prove strict convexity with respect to the Weyl-Petersson metric. In this section, we verify that this property continues to hold even when the curvature of $g$ is no-longer necessarily constant.
\par
Let $B$ be a trace-free Codazzi field of $h$ which we consider as a tangent vector of $\opTeich[S]$ at $h$ (c.f. \cite{Tromba}).
\proclaim{Lemma \nextprocno}
\noindent There exists $\epsilon>0$ and a smooth function $\phi:]-\epsilon,\epsilon[\times\Surface\rightarrow\Bbb{R}$ such that, for all $t\in]-\epsilon,\epsilon[$,
$$
h_t := h([(1+t^2\phi)\opId  + tB]\cdot,[(1+t^2\phi)\opId  + tB]\cdot)\eqnum{\nexteqnno[DefinitionOfGeodesic]}
$$
is a hyperbolic metric, where here $\phi_t:=\phi(t,\cdot)$. Furthermore,
$$
(\Delta - 2)\phi_0 = \opDet(B)\eqnum{\nexteqnno[FormulaForPhiZero]}
$$
and, up to and including order $2$ in $t$, $(h_t)_{t\in]-\epsilon,\epsilon[}$ defines a Weyl-Petersson geodesic in $\opTeich[S]$.
\endproclaim
\proclabel{WeylPeterssonGeodesic}
\proof For all $t\in\Bbb{R}$ and for all $\phi\in C^\infty(S)$, denote
$$\eqalign{
B_{t,\phi} &:= (1+t^2\phi)\opId + tB\ \text{and}\cr
h_{t,\phi} &:= h(B_{t,\phi}\cdot,B_{t,\phi}\cdot).\cr}
$$
Using \eqnref{DifferentialCurvatureFormulaA}, we show that, for all $(t,\phi)$ sufficiently close to $(0,0)$, the curvature of $h_{t,\phi}$ satisfies
$$
\kappa[h_{t,\phi}] = -1-t^2[(\Delta-2)\phi - \opDet(B)] + \opO(t^3).
$$
The first two assertions now follow by the implicit function theorem. To prove the final assertion, observe first that, since $B$ is trace-free and symmetric, its square is a scalar multiple of the identity, so that
$$
h_t = 2th(B\cdot,\cdot) + t^2(2\phi + \opTr(B^2)/2)h.
$$
On the other hand, since $B$ is a trace-free Codazzi field,
$$
h(B\cdot,\cdot) = Q + \overline{Q}
$$
for some quadratic, holomorphic differential $Q$. The metric $h_t$ thus has the form of $(5.1)$ of \cite{Wolf}. The final assertion now follows by Corollary $5.4$ of that paper.\qed
\medskip
We henceforth denote
$$
B_t := (\opId + t^2\phi_t) + t B.\eqnum{\nexteqnno[DefinitionOfBt]}
$$
For all $t$, let $\Phi_t:(\Sigma,h_t)\rightarrow(\Sigma,g)$ and $A_t$ be respectively the unique one-harmonic diffeomorphism and positive-definite Codazzi field of $g$ such that
$$
g_t := \Phi_t^*g = h_t(A_t\cdot,A_t\cdot).\eqnum{\nexteqnno[DefinitionOfgt]}
$$
By Theorem \procref{TrapaniValliII}, both $(A_t)_{t\in]-\epsilon,\epsilon[}$ and $(\Phi_t)_{t\in]-\epsilon,\epsilon[}$ vary smoothly with $t$. Finally, for all $t$, define
$$
\hat{\opE}_t := \hat{\opE}_g[h_t].\eqnum{\nexteqnno[FamilyOfHt]}
$$
In the present framework, Lemma $3.4$ of \cite{BonMonSch} becomes
\proclaim{Lemma \nextprocno}
\noindent Upon reducing $\epsilon$ if necessary, there exists a smooth family $(X_t)_{t\in]-\epsilon,\epsilon[}$ of vector fields and a smooth family $(f_t)_{t\in]-\epsilon,\epsilon[}$ of functions such that, for all $t$,
$$
\dot{A}_t + B_t^{-1}\dot{B}_tA_t = \nabla^tX_t + f_tJ^tA_t\eqnum{\nexteqnno[BasicRelation]}
$$
where, for all $t$, $\nabla^t$ and $J^t$ denote respectively the Levi-Civita covariant derivative and the complex structure of $h_t$.
\endproclaim
\proclabel{BasicRelation}
\noindent This yields
\proclaim{Lemma \nextprocno}
\noindent The first derivative of $\hat{E}_t$ satisfies
$$
\frac{\partial\hat{E}_t}{\partial t} = \int_S\opTr(A_t)\opTr(B_t^{-1}\dot{B}_t) - \opTr(A_tB_t^{-1}\dot{B}_t)\opdArea[h_t].\eqnum{\nexteqnno[FirstDerivativeOfF]}
$$
In particular, since $\dot{B}_0$ is trace-free,
$$
\frac{\partial\hat{E}_t}{\partial t}\bigg|_{t=0} = -\int_S\opTr(A_0\dot{B}_0)\opdArea[h_0].\eqnum{\nexteqnno[FirstDerivativeOfFAtZero]}
$$
\endproclaim
\proclabel{FirstDerivativeOfF}
\proof Indeed, by \eqnref{TwoDimMatrixRelations}, taking the trace of \eqnref{BasicRelation} yields
$$
\opTr(\dot{A}_t) + \opTr(A_tB^{-1}_t\dot{B}_t) = \nabla^{h_t}\cdot X_t,
$$
Using Stokes' Theorem, we therefore obtain
$$\eqalign{
\frac{\partial\hat{E}_t}{\partial t} &= \int_S\opTr(\dot{A}_t)\opdArea[h_t] + \int_S\opTr(A_t)\frac{\partial}{\partial t}\opdArea[h_t]\cr
&=-\int_S\opTr(A_tB^{-1}_t\dot{B}_t)\opdArea[h_t]+\int_\Surface\opTr(A_t)\opTr(B_t^{-1}\dot{B}_t)\opdArea[h_t],\cr}
$$
as desired.\qed
\medskip
\noindent The main difference between \cite{BonMonSch} and the present case is
\proclaim{Lemma \nextprocno}
\noindent Using the notation of Lemma \procref{BasicRelation},
$$
J^0X_0 = \opDet(A_0)A_0^{-1}\nabla^0f_0.\eqnum{\nexteqnno[UsefulRelationIII]}
$$
\endproclaim
\proclabel{UsefulRelationIII}
\remark[UsefulRelationIII] This result substitutes Lemma $3.8$ of \cite{BonMonSch}.
\medskip
\proof Since $A_t$ is an $h_t$-Codazzi field for all $t$,
$$
d^{\nabla^{h_t}}A_t = 0.\eqnum{\nexteqnno[ConvexityCodazziRelation]}
$$
Since $B_0=\opId$ and since $\dot{B}_0$ is an $h_0$-Codazzi field,
$$
\nabla^{h_t} = B_t^{-1}\nabla^{h_0}B_t + \opO(t^2),
$$
so that, at zero,
$$
\frac{\partial}{\partial t}\nabla^{h_t} = [\nabla^{h_0},\dot{B}_0].
$$
Differentiating \eqnref{ConvexityCodazziRelation} at zero therefore yields
$$
d^{\nabla^{h_0}}\big(A_0 + \dot{B}_0A_0\big) = \dot{B}_0d^{\nabla^{h_0}}A_0 = 0,
$$
so that, by \eqnref{BasicRelation},
$$
d^{\nabla^{h_0}}(\nabla^{h_0}X_0 + f_0 J^{h_0}A_0) = 0.
$$
Now let $(e_1,e_2)$ be an orthonormal frame of the metric $h_0$. Since $h_0$ is hyperbolic,
$$
d^{\nabla^{h_0}}(\nabla^{h_0}X_0)(e_1,e_2) = R^{h_0}_{e_1e_2}X_0 = J^{h_0}X_0.
$$
On the other hand
$$
d^{\nabla^{h_0}}(J^{h_0}A_0) = J^{h_0}d^{\nabla^{h_0}}A_0 = 0,
$$
so that, bearing in mind \eqnref{TwoDimMatrixRelations},
$$\eqalign{
d^{\nabla^{h_0}}(f_0J^{h_0}A_0)(e_1,e_2) &= (df_0\wedge(J^{h_0}A_0))(e_1,e_2)\cr
&= J^{h_0}A_0J^{h_0}\nabla^{h_0}f_0\cr
&= -\opDet(A_0)A_0^{-1}\nabla^{h_0}f_0.\cr}
$$
The result follows.\qed
\medskip
\noindent The remainder of Bonsante, Mondello \& Schlenker's argument may now be applied without further changes. We therefore obtain
\proclaim{Lemma \nextprocno}
\noindent The second derivative of $\opE_t$ satisfies
$$
\frac{\partial^2\hat{\opE}_t}{\partial t^2}\bigg|_0 \geq -\int_\Surface\opTr(J_0A_0J_0\ddot{B}_0)\opdArea[h_0].\eqnum{\nexteqnno[SecondDerivativeOfF]}
$$
\endproclaim
\proclabel{SecondDerivativeOfF}
\noindent This yields
\proclaim{Lemma \nextprocno}
\noindent $\hat{\opE}_g$ is strictly convex with respect to the Weyl-Petersson metric.
\endproclaim
\proclabel{StrictConvexityOfE}
\proof Indeed, by \eqnref{FormulaForPhiZero},
$$
\ddot{B}_0 = 2\phi_0\opId,
$$
where
$$
(\Delta-2)\phi_0 = \opDet(B).
$$
Since $B$ is symmetric and trace-free, its determinant is non-positive and vanishes if and only if $B$ vanishes. It follows that
$$
\Delta\phi_0 \leq 2\phi_0,
$$
so that, by the maximum principle, $\phi_0$ is non-negative. Furthermore, if $\phi_0$ vanishes then so too does $\opDet(B)$ and therefore also $B$. Since $B$ is non-zero, it follows that $\phi_0$ is strictly positive at at least one point, so that, by Lemma \procref{SecondDerivativeOfF},
$$
\ddot{\opE}_0 \geq -\int_\Surface\opTr(J_0A_0J_0\ddot{B}_0)\opdArea[h_0] \geq \int_\Surface\phi_0\opTr(A_0)\opdArea[h_0] > 0.
$$
Since $(h_t)_{t\in]-\epsilon,\epsilon[}$ is a Weyl-Petersson geodesic of $\opTeich[S]$ up to and including order $2$ in $t$, the result follows.\qed
\newsubhead{Proof of Theorem \procref{MainTheoremI}}[ProofOfMainTheorem]
We now prove our main result. Let $\Surface$ be a compact surface of hyperbolic type and let $g_\pm$ be negatively curved metrics over $\Surface$. Define the functional $\opF:\opTeich[\Surface]\rightarrow\Bbb{R}$ by
$$
\opF[h] := \hat{\opE}_{g_+}[h] + \hat{\opE}_{g_-}[h],\eqnum{\nexteqnno[SumOfConvexFunctionals]}
$$
where $\hat{\opE}$ is the functional defined in Section \subheadref{TheDerivativesOfE}. By Theorem \procref{CompactnessForCodazziFields} and Lemma \procref{StrictConvexityOfE}, $\opF$ is a proper function over $\opTeich(\Surface)$ which is strictly convex with respect to the Weyl-Petersson metric. Since the Weyl-Petersson metric is geodesically convex (c.f. Section $5.1$ of \cite{Wolf}), this function has a unique minimum at some point $h_0$, say. It remains to show that $h_0$ is the desired marked hyperbolic metric.
\par
Let $A_\pm$ denote the unique positive-definite, symmetric Codazzi field of $h_0$ such that $h_0(A_\pm\cdot,A_\pm\cdot)=\Phi_\pm^*g_\pm$ for some one-harmonic diffeomorphism $\Phi_\pm$ which preserves the marking. Since $h$ is a critical point of $F$, by \eqnref{FirstDerivativeOfFAtZero},
$$
\int_\Surface\opTr((A_++A_-)B)\opdArea = 0,\eqnum{\nexteqnno[SumOfAOneAndATwoIsExact]}.
$$
for every trace-free Codazzi field $B$ of $h_0$. It follows by Lemma $3.3$ of \cite{FillSmiII} that
$$
A_+ + A_- = f\opId - \opHess(f),\eqnum{\nexteqnno[PrimitiveOfSumOfAOneAndATwo]}
$$
for some smooth function $f:\Surface\rightarrow\Bbb{R}$, where here $\opHess$ denotes the Hessian operator of $h_0$.
\par
Let $\rho:\pi_1(\Surface)\rightarrow\opPSL(2,\Bbb{R})$ be such that $(\Surface,h_0)$ identifies with $\Bbb{H}^2/\rho(\pi_1(\Surface))$ as a point of Teichm\"uller space. We identify $\Bbb{H}^2$ with the future component of the unit pseudo-sphere in $\Bbb{R}^{2,1}$. Likewise, we identify $T\Bbb{H}^2$ with a subbundle of the trivial bundle $\Bbb{H}^2\times\Bbb{R}^{2,1}$. Finally, we identify $f$ with its lift to a function over $\Bbb{H}^2$, and we indentify $A_\pm$ with its lift to a section of $\opEnd(T\Bbb{H}^2)$. In particular, $A_\pm$ identifies with a section of  $T^*\Bbb{H}^2\otimes\Bbb{R}^{2,1}$. Given a base point $x_0\in\Bbb{H}^2$, we now define $X_\pm:\Bbb{H}^2\rightarrow\Bbb{R}^{2,1}$ by
$$
X_\pm := U_\pm \pm \int_{x_0}^x(A_\pm\cdot\partial_\tau)d\tau,\eqnum{\nexteqnno[DefinitionOfXOneAndXTwo]}
$$
where the vectors $U_\pm$ are chosen such that
$$
U_+ - U_- = \nabla f(x_0).\eqnum{\nexteqnno[CharacteristicPropertyOfU]}
$$
Indeed, since $A_\pm$ is a Codazzi tensor, the integrals \eqnref{DefinitionOfXOneAndXTwo} is well-defined independent of the path chosen from $x_0$ to $x$.
\par
The functions $X_+$ and $X_-$ define respectively future- and past-oriented, strictly convex embeddings of $\Bbb{H}^2$ into $\Bbb{R}^{2,1}$. The metric induced over $\Bbb{H}^2$ by $X_\pm$ is trivially equal to $h_0(A_\pm\cdot,A_\pm\cdot)$. In addition, $X_\pm$ is equivariant with respect to the affine deformation $(\rho,\tau_\pm)$, where
$$
\tau_\pm := \pm\int_{x_0}^{\gamma(x_0)}(A_\pm\cdot\partial_\tau)d\tau \pm U_1\mp\rho(\gamma)(U_\pm).
\eqnum{\nexteqnno[DefinitionOfRhoOneAndRhoTwo]}
$$
In order to prove existence it thus suffices to show that $\tau_+=\tau_-=:\tau$. However, the support function of $X_\pm$ is
$$
\phi_\pm(x) := \pm\langle X_\pm(x),x\rangle.\eqnum{\nexteqnno[DefinitionOfSupportFunctions]}
$$
It follows upon integrating \eqnref{PrimitiveOfSumOfAOneAndATwo} that
$$
\phi_+ + \phi_- = f,
$$
so that, since $f$ is bounded,
$$
\mlim_{x\rightarrow\partial_\infty\Bbb{H}^2}\frac{\phi_+(x)}{\opCosh(d(x,x_0))} = \mlim_{x\rightarrow\partial_\infty\Bbb{H}^2}\frac{\phi_-(x)}{\opCosh(d(x,x_0))},\eqnum{\nexteqnno[LimitsOfPhisCoincide]}
$$
where $d(x,x_0)$ here denotes the distance in $\Bbb{H}^2$ from $x_0$ to $x$. Recall now that $\phi_\pm$ defines a continuous function $\tilde{\phi}_\pm$ over $\partial_\infty\Bbb{H}^2$ (c.f. \cite{BarbFill}). By \eqnref{LimitsOfPhisCoincide}, $\tilde{\phi}_+=\tilde{\phi}_-$ so that, by Corollary $3.14$ of \cite{BarbFill}, $\tau_+=\tau-$ and existence follows.
\par
Finally, let $(\rho,\tau)\in\opSO_0(2,1)\ltimes\Bbb{R}^{2,1}$ be an affine deformation, let $e_\pm:\tilde{S}\rightarrow\Bbb{R}^{2,1}$ be LSC, $(\rho,\tau)$-equivariant, spacelike immersions such that $e_+$ is future-oriented, $e_-$ is past-oriented, and
$$
g_{\pm} = e_{\pm}^*\delta^{2,1}.
$$
Let $A_\pm$ be the shape operator of $e_\pm$ and recall that $A_\pm$ is a Codazzi field. By Lemma $3.3$ of \cite{FillSmiII}, there exist smooth functions $f_\pm:S\rightarrow\Bbb{R}$ and trace-free Codazzi fields $B_\pm$ such that
$$
A_\pm = B_\pm + (f_\pm\opId - \opHess_h(f_\pm)).
$$
Observe that the cocycle of $A_\pm$ is $\pm\tau$ so that the cocycle of $A_++A_-$ and therefore also of $B_++B_-$ vanishes. It follows by Proposition $3.17$ of \cite{BarbFill} (c.f. also \cite{BonSeppi}) that $B_++B_-$ vanishes. Consequently, for every trace-free Codazzi field $B$,
$$
\int_S\opTr((A_++A_-)B)\opdArea = 0,
$$
so that $h$ is the unique critical point of $F$, and uniqueness follows.
\inappendicestrue
\global\headno=0
\bigskip
\goodbreak
\newhead{The generalized lorentzian metric}[TheGeneralisedLorentzianMetric]
\noindent In this appendix, we show how the generalised lorentzian metric studied by Trapani \& Valli in Section $5$ of \cite{TrapaniValli} has a simpler expression in the terminology of the present paper. First, let $b$ be the non-degenerate, symmetric bilinear form defined over $\opEnd(2)$ by
$$
b(B,B) := \frac{1}{2}\opTr(BJB^tJ) = -\opDet(B),\eqnum{\nexteqnno[MinkowskiBilinearForm]}
$$
where $J$ here denotes the standard complex structure of $\Bbb{R}^2$. Since $b$ is invariant under the conjugation action of $\opGL(2)$ on $\opEnd(2)$, it extends to a bi-invariant, semi-riemannian metric over $\opGL(2)$ which we also denote by $b$.
\par
The orthogonal complement of $\frak{so}(2)$ in $\opEnd(2)$ with respect to $b$ is $\opSymm(2)$. Since
$$\eqalign{
[\frak{so}(2),\opSymm(2)] &\subseteq \opSymm(2)\ \text{and}\cr
[\opSymm(2),\opSymm(2)] &\subseteq \frak{so}(2),\cr}
$$
the quotient $\opSymm_+(2)=\opGL(2)/\opSO(2)$ is a symmetric space. Furthermore, $b$ descends to a lorentzian metric over this space which we also denote by $b$. Identifying $\opSymm_+(2)$ with the space of positive-definite, symmetric matrices in the natural manner, we obtain
$$
b(A)(B,B) = \frac{1}{8}\opDet(A)^{-1}\opTr(BJBJ) = -\frac{1}{4}\opDet(A)^{-1}\opDet(B).\eqnum{\nexteqnno[LorentzMetricInQuotient]}
$$
Let $\Omega:=\nabla-D$ denote the difference between the Levi-Civita covariant derivative of $b$ and the standard derivative of $\opEnd(2)$. By the Koszul formula, when $[A,B]=0$,
$$
\Omega(A)(B,B) = -A^{-1}B^2.\eqnum{\nexteqnno[MetricOverQuotient]}
$$
\par
For $\alpha\in\Bbb{R}$, let $b_\alpha$ be the conformally rescaled metric defined over $\opSymm_+(2)$ by
$$
b_\alpha(A)(B,B) := \opDet(A)^\alpha b,\eqnum{\nexteqnno[ConformallyRescaledMetric]}
$$
and let $\Omega_\alpha := \nabla^\alpha - D$ denote the difference between its Levi-Civita covariant derivative and the standard derivative of $\opEnd(2)$. By the Koszul formula again, for all $\alpha$, when $[A,B]=0$,
$$
\Omega_\alpha(A)(B,B) = (\alpha-1)A^{-1}B^2.\eqnum{\nexteqnno[CovDerOfRescaledMetric]}
$$
In particular, when $\alpha=1/2$, the geodesics of $b_{1/2}$ passing through $\opId$ are all of the form
$$
\gamma(t) = (\opId + tA)^2,\eqnum{\nexteqnno[GeodesicsOfRescaledMetric]}
$$
for some $A\in\opSymm(2)$.
\par
Let $\Omega$ be an open subset of $\Bbb{R}^2$, let $\opSymm_+(2,\Omega)$ be as in Section \subheadref{TheEnergyDensity}, and let $g_0\in\opSymm_+(2,\Omega)$ be a fixed metric over $\Omega$. As in Section $5$ of \cite{TrapaniValli}, let $\langle\cdot,\cdot\rangle$ be the semi-riemannian metric defined over $\opSymm_+(2,\Omega)$ such that, for all $g\in\opSymm_+(2,\Omega)$ and for all $h\in\opSymm(2,\Omega)$,
$$
\langle h,h\rangle_g := -\frac{1}{2}\int_\Omega\opDet(g^{-1}h)\opdVol_g = -\frac{1}{4}\int_\Omega\opDet(g_0g^{-1})^{-1/2}\opDet(g_0^{-1}g)\opdVol_{g_0}.\eqnum{\nexteqnno[DefnOfGeneralisedLorentzMetric]}
$$
Observe that the geodesic equation for $\langle\cdot,\cdot\rangle$ may be solved pointwise so that, by \eqnref{GeodesicsOfRescaledMetric}, its geodesics passing through $g_0$ are precisely those curves of the form
$$
g_t := g_0((\opId + tA)^2\cdot,\cdot) = g_0((\opId + tA)\cdot,(\opId + tA)\cdot),\eqnum{\nexteqnno[GeodesicOfLorentzianMetric]}
$$
where $A$ is an endomorphism field over $\Omega$ which is symmetric with respect to $g_0$. The exponential map about $g_0$ of the lorentzian metric is thus
$$
\Psi:U\rightarrow\opSymm_+(2,\Omega);A\mapsto g_0((\opId + tA)\cdot,(\opId + tA)\cdot),\eqnum{\nexteqnno[ExponentialMap]}
$$
where
$$
U := \left\{ A\ |\ A(x) + \opId \geq 0\ \forall x\right\}.\eqnum{\nexteqnno[DomainOfExponentialMap]}
$$
\par
It remains to compare this construction with that of Section $5$ of \cite{TrapaniValli}. We may assume that $g_0$ is conformal to the standard metric over $\Bbb{R}^2$. That is,
$$
g_0 = \opRe(\rho^2dzd\overline{z}) = \frac{1}{2}\rho^2(dzd\overline{z} + d\overline{z}dz),
$$
for some positive function $\rho$. Given a real number $P$ and a Beltrami differential $Q$ of the form
$$
Q = (a+bi)\frac{d\overline{z}}{dz} = (a+ bi)\partial_z\otimes d\overline{z},
$$
let $\tilde{A}(P,Q)$ be the matrix defined by
$$
\tilde{A}(P,Q) := Q + P\opId + \overline{Q} = \pmatrix a+P\hfill& b\hfill\cr b\hfill& -a+P\hfill\cr\endpmatrix.\eqnum{\nexteqnno[DefinitionOfATilde]}
$$
Trivially, $\tilde{A}$ depends linearly on $P$ and $Q$. Furthermore
$$\eqalign{
\opTr(\opId+\tilde{A}(P,Q)) &= 2(P+1)\ \text{and}\hfill\cr
\opDet(\opId + \tilde{A}(P,Q)) &= (P+1)^2 - \left|Q\right|^2.\hfill\cr}
$$
The preimage of $U$ under $\tilde{A}$ is thus
$$
\tilde{U} := \left\{ (P,Q)\ |\ P+1> 0\ \text{and}\ (P+1)^2 - \left|Q\right|^2 > 0\right\},\eqnum{\nexteqnno[DefinitionOfTildeU]}
$$
which is precisely the set denoted by $\Omega$ in Section $5$ of \cite{TrapaniValli}. The composition $\tilde{\Psi}:=\Psi\circ\tilde{A}$ is
$$
\tilde{\Psi}(P,Q) := \big((P+1)Q + ((P+1)^2 + \left|Q\right|^2) + (P+1)\overline{Q}\big)g,\eqnum{\nexteqnno[DefinitionOfTildePsi]}
$$
which is precisely the function denoted by $\Psi$ in Section $5$ of \cite{TrapaniValli}.
\goodbreak
\newhead{Bibliography}[Bibliography]
{\leftskip = 5ex \parindent = -5ex
\leavevmode\hbox to 4ex{\hfil \cite{Ahlfors}}\hskip 1ex{Ahlfors L., {\sl Lectures on quasiconformal mappings}, AMS, (2006)}
\medskip
\leavevmode\hbox to 4ex{\hfil \cite{AndEtAl}}\hskip 1ex{Andersson L., Barbot T., Benedetti R., Bonsante F., Goldman W. M., Labourie F., Scannell K. P., Schlenker J. M., Notes on: “Lorentz spacetimes of constant curvature”, {\sl Geom. Dedicata}, {\bf 126}, (2007), 47--70}
\medskip
\leavevmode\hbox to 4ex{\hfil \cite{BarbFill}}\hskip 1ex{Barbot T., Fillastre F., Quasi-Fuchsian co-Minkowski manifolds, to appear {\sl In the tradition of Thurston}, Springer-Verlag, (2020)}
\medskip
\leavevmode\hbox to 4ex{\hfil \cite{Bers}}\hskip 1ex{Bers L., Simultaneous uniformization, {\sl Bull. Amer. Math. Soc.}, {\bf 66}, no.2, (1960), 94--97}
\medskip
\leavevmode\hbox to 4ex{\hfil \cite{BonMonSch}}\hskip 1ex{Bonsante F., Mondello G., Schlenker J. M., A cyclic extension of the earthquake flow II, {\sl Ann. Sci. Ec. Norm. Sup\'er.}, {\bf 48}, (2015), no. 4, 811--859}
\medskip
\leavevmode\hbox to 4ex{\hfil \cite{BonSeppi}}\hskip 1ex{Bonsante F., Seppi A., On Codazzi tensors on a hyperbolic surface and flat Lorentzian geometry, {\sl Int. Math. Res. Not. IMRN}, (2):343–417, (2016)}
\medskip
\leavevmode\hbox to 4ex{\hfil \cite{Buser}}\hskip 1ex{Buser P., {\sl Geometry and spectra of compact Riemann surfaces}, Birkh\"auser Verlag, (1992)}
\medskip
\leavevmode\hbox to 4ex{\hfil \cite{FillSmi}}\hskip 1ex{Fillastre F., Smith G., Group actions and scattering problems in Teichm\"uller theory, in {\sl Handbook of group actions IV}, Advanced Lectures in Mathematics, {\bf 40}, (2018), 359--417}
\medskip
\leavevmode\hbox to 4ex{\hfil \cite{FillSmiII}}\hskip 1ex{Fillastre F., Smith G., A note on invariant constant curvature immersions in Mink\-ow\-ski space, to appear in {\sl Geom. Dedicata}}
\medskip
\leavevmode\hbox to 4ex{\hfil \cite{GilbTrud}}\hskip 1ex{Gilbarg D., Trudinger N. S., {\sl Elliptic partial differential equations of second order}, Classics in Mathematics, Springer-Ver\-lag, Berlin, (2001)}
\medskip
\leavevmode\hbox to 4ex{\hfil \cite{Jost}}\hskip 1ex{Jost J., {\sl Riemannian Geometry and Geometric Analysis}, Universitext, Springer-Verlag, (2011)}
\medskip
\leavevmode\hbox to 4ex{\hfil \cite{LabI}}\hskip 1ex{Labourie F., Probl\`eme de Minkowski et surfaces \`a courbure constante dans les vari\'et\'es hyperboliques, {\sl Bull. Soc. math. France}, {\bf 119}, (1991), 307--325}
\medskip
\leavevmode\hbox to 4ex{\hfil \cite{LabII}}\hskip 1ex{Labourie, F., Metriques prescrites sur le bord des vari\'et\'es hyperboliques de dimension $3$, {\sl J. Differ. Geom.}, {\bf 35}, (1992), 609--626}
\medskip
\leavevmode\hbox to 4ex{\hfil \cite{Mess}}\hskip 1ex{Mess G., Lorentz spacetimes of constant curvature, {\sl Geom. Dedicata}, {\bf 126}, (2007), 3--45}
\medskip
\leavevmode\hbox to 4ex{\hfil \cite{Mumford}}\hskip 1ex{Mumford D., A remark on Mahler's compactness theorem, {\sl Proc. AMS}, {\bf 28}, no.1, (1971), 289--294}
\medskip
\leavevmode\hbox to 4ex{\hfil \cite{NewNir}}\hskip 1ex{Newlander A., Nirenberg L., Complex analytic coordinates in almost complex manifolds, {\sl Ann. of Math.}, {\bf 65}, (1957), 391--404}
\medskip
\leavevmode\hbox to 4ex{\hfil \cite{Schlenker}}\hskip 1ex{Schlenker J. M., Hyperbolic manifolds with convex boundary, {\sl Inventiones mathematicae}, {\bf 163}, (2006), 109--169}
\medskip
\leavevmode\hbox to 4ex{\hfil \cite{Schoen}}\hskip 1ex{Schoen R., The role of harmonic mappings in rigidity and deformation problems, in {\sl Complex geometry: Proceedings of the Osaka international conference}, (1993), 179-200}
\medskip
\leavevmode\hbox to 4ex{\hfil \cite{SchoenYau}}\hskip 1ex{Schoen R., Yau S. T., On univalent harmonic maps between surfaces, {\sl Invent. Math.}, {\bf 44}, (1978), 265--278}
\medskip
\leavevmode\hbox to 4ex{\hfil \cite{Tamburelli}}\hskip 1ex{Tamburelli A., Prescribing metrics on the boundary of anti-de Sitter $3$-manifolds, {\sl Int. Math. Res. Not. IMRN}, (2018), no. 5, 1281--1313}
\medskip
\leavevmode\hbox to 4ex{\hfil \cite{Thurston}}\hskip 1ex{Thurston W., The geometry and topology of three manifolds}
\medskip
\leavevmode\hbox to 4ex{\hfil \cite{TrapaniValli}}\hskip 1ex{Trapani S., Valli G., One-harmonic maps on Riemann surfaces, {\sl Comm. Anal. Geom.}, {\bf 3}, no. 4, (1985), 645--681}
\medskip
\leavevmode\hbox to 4ex{\hfil \cite{Tromba}}\hskip 1ex{Tromba A. J., {\sl Teichm\"uller theory in Riemannian geometry}, Lectures in Mathematics, ETH Z\"urich, Birkh\"auser Verlag, Basel, (1992)}
\medskip
\leavevmode\hbox to 4ex{\hfil \cite{Wolf}}\hskip 1ex{Wolf M., The Teichm\"uller theory of harmonic maps, {\sl J. Diff. Geom}, {\bf 29}, (1989), 449--479}
\par}
%
%
%
%
\enddocument